\documentclass[12pt]{amsart}

\usepackage{mathrsfs}
\usepackage{amsmath,amssymb,enumerate,tikz}
\usepackage{amsfonts,amssymb}
\usepackage{booktabs}
\usepackage[colorlinks,linkcolor=red,anchorcolor=blue,
citecolor=green,CJKbookmarks=True]{hyperref}
\usepackage{xypic}
\usepackage[]{mdframed}
\usepackage[T1]{fontenc}
\usepackage{csquotes}
\usepackage{soul, xcolor, todonotes, easyReview}

\theoremstyle{plain}
\newtheorem{theorem}{Theorem}[section]
\newtheorem{lemma}[theorem]{Lemma}
\newtheorem{corollary}[theorem]{Corollary}

\newtheorem{proposition}[theorem]{Proposition}

\theoremstyle{definition}
\newtheorem{remark}[theorem]{Remark}
\newtheorem{question}[theorem]{Question}

\newtheorem{convention}{Convention}
\newtheorem{example}{Example}

\newcounter{consta}
\renewcommand{\theconsta}{{\kappa_{\arabic{consta}}}}
\newcounter{constb}[section]

\newcounter{constc}
\renewcommand{\theconstc}{{c_{\arabic{constc}}}}
\newcounter{constA}

\newcommand{\consta}{\refstepcounter{consta}\theconsta}

\newcommand{\constc}{\refstepcounter{constc}\theconstc}

\newcommand{\tec}{Teichm\"uller }

\newcommand{\C}{\mathscr{C}}
\newcommand{\T}{\mathcal{T}}
\newcommand{\M}{\mathcal{M}}

\newcommand{\R}{\mathbb R}
\newcommand{\Z}{\mathbb{Z}}
\renewcommand{\H}{\mathbb{H}}

\newcommand{\Int}{\mathrm{Int}}
\newcommand{\Vol }{\mathrm{Vol}}
\renewcommand{\i }{\mathrm{i}}
\newcommand{\sys}{\mathrm{sys}}
\newcommand{\thin }{\mathrm{thin}}
\newcommand{\thick  }{\mathrm{thick}}

\newcommand{\CH}{\mathrm{CH}}
\newcommand{\arcsinh}{\mathrm{arcsinh}~}

\title[Algebraic intersection for hyperbolic surfaces]{Algebraic intersection for hyperbolic surfaces}

\author{Manman Jiang \and Huiping Pan}
\date{\today}
\address{ Manman Jiang\\
 Guangzhou Maritime University,  Guangzhou, 510725, China}
\email{chnjiangmm@foxmail.com}

\address{Huiping Pan\\ School of mathematics, South China University of Technology, 510641, Guangzhou, China}
\email{panhp@scut.edu.cn}

\begin{document}

\maketitle
\begin{abstract}
We show that the algebraic intersection form of hyperbolic surfaces of genus $g$ has a minimum in the moduli space and that the minimum grows in the order $(\log g)^{-2}$  in terms of the genus.  We also describe the asymptotic behavior of the algebraic intersection form in the moduli space as the homologically systolic length goes to zero.
\medskip

\end{abstract}

\section{Introduction}

\subsection{Main results}
Given a Riemannian surface $X$, the \emph{algebraic intersection form} $K(X)$ is defined to be: 
\begin{equation}\label{eq:intection:length}
	K(X)=\sup_{\alpha,\beta}\frac{|\Int(\alpha,\beta)|}{\ell_X(\alpha)\ell_X(\beta)},
\end{equation}
where $\Int(\cdot,\cdot)$ represents the \emph{algebraic intersection number},
where $\alpha$ and $\beta$ range over all oriented simple closed geodesics under the Riemannian metric $X$, and where $\ell_X(\cdot )$ represents the length of geodesics with respect to $X$. The algebraic intersection form $K(X)$ and its volume $\Vol(X)$ of a Riemannian surface $X$ satisfies $K(X)\Vol(X)\geq 1$ with equality holds if and only if $X$ is a flat torus \cite{MM2014}. 
     
     For hyperbolic surfaces, Massart and Muetzel \cite{MM2014} proved that $K(X)$ is unbounded and nonproper in the moduli space $\M_g$ of hyperbolic surfaces of genus $g\geq2$.  In particular, $K(X)$ does not have a maximum in $\M_g$. Regarding the minimum of $K(X)$, they asked:
     \begin{question}[\cite{MM2014}, Question 1.7]
     	Does $K(X)$ have a minimum in the moduli space $\M_g$? If so, which surface realizes the minimum? 
     \end{question}
   
In this paper, we answer the first part of the above question affirmatively:
\begin{theorem}\label{thm:minimum}For $g\geq2$, the algebraic intersection form $K(X)$ attains its   minimum in the moduli space $\M_g$. 
\end{theorem}

Next, we estimate the minimum in terms of the genus:
\begin{theorem} \label{thm:min:estimate} Let $g\geq2$. 
There exist  two universal constants $C_1, C_2$ such that 
	\begin{equation*}
	\frac{C_1}{(\log g)^2}\leq \min_{X\in\M_g}K(X)\leq \frac{C_2}{(\log g)^2}.
	\end{equation*} 
\end{theorem}
\begin{remark}
In \cite[Section 1]{Boulanger2023b}, Boulanger also discussed the minimum, or more precisely the inifimum, and sketched a proof of the upper bound. The upper bound also follows directly from \cite[Proposition 1.2]{Torkaman2023}. Our main contribution here is the lower bound.
\end{remark}

This leads to the question:
\begin{question}
	Does the limit $\lim\limits_{g\to\infty} (\log g)^2 \min\limits_{X\in\M_g} K(X)$ exist as $g\to\infty$?
\end{question}
	
	Another interesting question is to study the asymptotic behavior in the moduli space. Massart and Muetzel \cite{MM2014} observed that it is possible to pinch a homologically trivial simple closed geodesic while keeping the algebraic intersection form bounded from above. On the other hand, pinching a homologically nontrivial simple closed geodesic would alway make the algebraic intersection form diverge to infinity \cite{MM2014}:
	\begin{theorem}[\cite{MM2014}, Theorem 1.6]\label{thm:MM2014}
		There exist constants $A$ and $B$ depending only on the genus $g$ such that for any $X\in \M_g$ with the homologically systolic length $\sys_h(X)$ sufficiently small,  we have
		\begin{equation*}
			\frac{A}{\sys_h(X)\cdot |\log \sys_h(X)|}\leq K(X)\leq \frac{B}{\sys_h(X)\cdot |\log \sys_h(X)|}.
		\end{equation*}
	\end{theorem}

	Here the \emph{homologically systolic length} $\sys_h(X)$ is the length of any shortest homologically nontrivial simple closed geodesic on $X$. Based on Theorem \ref{thm:MM2014}, it is natural to wonder if it is possible to describe the asymptotic behavior in terms of homologically systolic length. However, this will not work. For any $g\geq2$ and for any $1\leq c\leq 2$,  we construct a sequence of hyperbolic surfaces $\{X_n\}$ in $\M_g$ such that, as $n\to\infty$, we have 
	\begin{equation*}
		K(X_n)\sim \frac{1}{2c\cdot \sys_h(X_n)|\log \sys_h(X_n)|}.
	\end{equation*}
	 Here $A\sim B$ means $\lim A/B=1$. See Section \ref{sec:exapmple} for the construction.  This example demonstrates that the asymptotic behavior of the algebraic intersection form cannot be described as a single function of the homologically systolic length. Instead, all short curves should be taken into consideration. 

\begin{theorem}\label{thm:asymp}
 For $g\geq2$ and $X\in\M_g$, let $\Gamma$ be the set of nonseparating simple closed geodesics of length at most $1$. For each $\gamma\in\Gamma$, let $\alpha_\gamma$ be one of the shortest nonseparating simple closed geodesics that intersect $\gamma$ exactly once. Then  as $\sys_h(X)\to0$, 
      \begin{equation*}
      K(X)\sim \max_{\gamma\in\Gamma} \frac{1}{\ell_X(\gamma)\ell_X(\alpha_\gamma)}.
      \end{equation*}
		\end{theorem}
		Here the constant $1$ in Theorem \ref{thm:asymp} can be replaced by any positive constant. The set $\Gamma$ may contain curves that are disjoint from any homological systole and cannot be replaced by the set of curves intersecting homological systoles (see Example \ref{example:disjoint:systole} in Section \ref{sec:exapmple}).

\subsection{Background and motivation}
Motivated by the question of how much can two curves of a given length intersect, Massart and Muetzel \cite{MM2014} initiated the study of \emph{algebraic intersection form} of surfaces with Riemannian metrics.  The algebraic intersection form plays an important role in the comparison between the stable norm and the Hodge norm ($L^2$-norm) \cite[Theorem 1.1]{MM2014}.  In the same paper, Massart and Muetzel also estimated the algebraic intersection form in terms of the homologically systolic length and the diameter.

	 Recently, the algebraic intersection form of translation surfaces has attracted much attention.  Cheboui,  Kessi, and Massart \cite{CKM2021a, CKM2021b} considered the scaling invariant quantity $\mathrm{KVol}(X):=K(X)\Vol(X)$ in the stratum $\mathcal{H}(2)$ of translation surfaces of genus two with a single cone point.  They computed the exact values of $\mathrm{KVol}$ for some family of square-tiled surfaces, proved that  $\inf_{X\in\mathcal{H}(2)}\mathrm{KVol}(X)\leq 2$ and constructed a sequence of square-tiled surfaces  $X_n$ in $\mathcal{H}(2)$  such that $\mathrm{KVol}(X_n)\to 2$ with 2 being the conjectural infimum for $\mathrm{KVol}$ over $\mathcal{H}(2)$. Boulanger \cite{Boulanger2023b} generalized this picture to the minimal stratum $\mathcal{H}(2g-2)$ of genus $g$ with a single cone point, proving that the infimum of $\mathrm{KVol}$ over each component of $\mathcal{H}(2g-2)$ is at most $g$ and constructing surfaces with $\mathrm{KVol}$ converging to $g$ with $g$ being the conjectural infimum. Beyond square-tiled surfaces, Boulanger, Lanner and Massart computed exact values of $\mathrm{KVol}$ for some non-arithmetic Veech surfaces,  the regular double $n$-gon translation surfaces for $n\geq5$ odd \cite{BLM2022} and regular $n$-gon translation surfaces for $n\geq8$ even \cite{Boulanger2023a}.

Using \emph{geometric} intersection number $\i(\cdot,\cdot)$, Torkaman \cite{Torkaman2023} initiated the study of the \emph{intersection strength} $I(X)$,  defined by 
 \begin{equation*}
 	I(X)=\sup_{\alpha,\beta}\frac{\i(\alpha,\beta)}{\ell_X(\alpha)\ell_X(\beta)}
 \end{equation*} 
 where the supremum ranges over all \emph{closed} geodesics on $X$. She proved that $I(X)$ is a proper function on the moduli space $\M_g$,  described the asymptotic behavior in terms of the systolic length $\sys(X)$ as $X$ approaches the boundary of $\M_g$,
 \begin{equation*}
 	I(X)\sim \frac{1}{2\cdot \sys(X)|\log \sys(X)|},
 \end{equation*}
  and showed that the minimum of $I(X)$ in the moduli space $\M_g$ also grows in the order $(\log g)^{-2}$ (by showing that every hyperbolic surface of genus $g$ has a figure eight closed geodesic of length at most $c\log g$ for some universal constant $c$).  From the definition of $K(X)$ and $I(X)$, we see that $K(X)\leq I(X)$, and hence $\min\limits_{X\in\M_g}K(X)\leq \min\limits_{X\in\M_g} I(X)$. It is natural to ask whether they are equal:
 \begin{question}
 	Are $\min\limits_{X\in \M_g} K(X)$ and $\min\limits_{X\in \M_g} I(X)$ equal to each other? 
 \end{question}

Another closely related quantity is the \emph{systolic volume}/\emph{area} $\mathrm{SysVol}(X)$ defined by:
 \begin{equation*}
 	\mathrm{SysVol}(X):=\frac{\Vol(X)}{\sys(X)^2},
 \end{equation*} 
 where $\sys(X)$ is the systolic length of $X$, that is, the length of any homotopically nontrivial shortest simple closed geodesic on $X$. For recent results about this quantity, we refer to  \cite{KS2005,KS2006,LiSu,KS2024}. 
  
  \subsection*{Sketch of the proof} The possible obstruction for the existence of a minimum of the algebraic intersection form in the moduli space is the existence of surfaces with short and homologically trivial simple closed geodesics.  We begin with an observation that, for the surfaces mentioned above, the algebraic intersection form is the maximum among that of the complementary components of short curves (Lemma \ref{lem:near:bnd}). We then apply the \emph{strip deformation} to lengthen those short curves and simple closed geodesics disjoint from them (Lemma \ref{lem:decreasing}), pushing the original surface into a compact subset of the moduli space and decreasing the algebraic intersection form. This ensures the existence of a minimum.  
  
     As mentioned earlier, the upper bound of Theorem \ref{thm:min:estimate} is already known in literature \cite{Boulanger2023b,Torkaman2023}. For the lower bound, the estimate relies crucially on the existence of short homological bases \cite{BPS2012} if the homological systole of the underlying surface is bounded from below. This allows us to reduce the problem to considering surfaces with \emph{short but not too short} homological systoles. We then apply the \emph{Fenchel-Nielsen coordinates} to lengthen those short curves to some uniform constant while keeping other simple closed geodesics from being too contracted (Lemma \ref{lem:auxiliary:surface}), yielding the desired lower bound.   
     
     For the asymptotic behavior in the moduli space, the key is to find good representatives within the homology classes of simple closed geodesics:  crossing each short curve at most once (Lemma \ref{lem:intersection:11}) and having small twists around each short curve ( Lemma \ref{lem:small:twist:11}).
     
      \subsection*{Acknowledgement} The first author is supported by the National Natural Science Foundation of China NSFC 12371076. The Second author is supported by National Natural Science Foundation of China NSFC 12371073 and Guangzhou Basic and Applied Basic Research Foundation 2024A04J3636.
 
 \tableofcontents

\section{Strip deformations}
\subsection{Teichm\"uller space}
Let $S$ be a closed orientable surface of genus $g\geq2$. The \tec space $\T_g$ of $S$ is the equivalence classes of hyperbolic metrics on $S$, where two metrics are said to be \emph{equivalent} if there exists an isometry between them which is homotopic to the identity map. The \tec space is topologically homeomorphic to $\R^{6g-6}$. The moduli space $\M_g$ is also defined as the equivalence classes of hyperbolic metrics on $S$, where two metrics are said to be \emph{equivalent} if there exists an isometry between them. The moduli space is the quotient space of the \tec space under the mapping class group.

 Given a hyperbolic surface $X$ and a simple closed curve $\alpha$, we denote by $\ell_X(\alpha)$ the length of the unique geodesic homotopic to $\alpha$. The shortest homotopically nontrivial closed geodesics on $X$ are called \emph{systoles} of $X$, and the shortest homologically nontrivial closed geodesics are called \emph{homological systoles} of $X$. We denote by $\sys(X)$ and $\sys_h(X)$  respectively the lengths of systoles and homologically systoles of $X$. 
     
The  \emph{geometric intersection number} $\i(\alpha,\beta)$ of two closed curves $\alpha$ and $\beta$ on $S$ is defined to be  the minimum of the cardinality of $\alpha'\cap\beta'$ where $\alpha'$ and $\beta'$ range over all curves homotopic to $\alpha$ and $\beta$, respectively.
 
 Given   a volume form $\omega$ on $S$, he  \emph{algebraic intersection number} $\Int(\alpha,\beta)$ of two differentiable oriented closed curves $\alpha:[0,1]\to S$ and $\beta:[0,1]\to S$ is defined as follows:
 \begin{equation*}
 	\Int(\alpha,\beta):=\sum \frac{\omega(\dot{\alpha}(s),\dot{\beta}(t))}{|\omega(\dot{\alpha}(s),\dot{\beta}(t))|},
 \end{equation*} 
 where $\dot{\alpha}$ and $\dot{\beta}$ represent the tangent vector of $\alpha$ and $\beta$ respectively and where the sum ranges over all pairs $(s,t)$ such that $\alpha(s)=\beta(t)$.
 
For more information about hyperbolic surfaces, Teichm\"uller space, and moduli space, we refer to \cite{IT1992,FM2012,Buser2010}.

\subsection{Strip deformations}
In this section, we introduce the
 \emph{strip deformations}, first introduced by Thurston (\cite[proof of Lemma 3.4]{Thurston1998}, see also \cite[Section 1.2]{DGK2016} and \cite{PapadopoulosTheret2010}) 

A \emph{strip} of width $w\geq0$ is a region on $\H^2$ bounded by two complete geodesic lines whose distance is exactly $w$. The two closest points on distinct components of the strip are called the \emph{feet} of the strip. 

Let $Y$ be a complete hyperbolic surface of finite genus and infinite area whose convex hull $\CH(Y)$ is a boarded hyperbolic surface. Let $\underline{\eta}=\{\eta_i\}$ be a collection of pairwise disjoint simple, proper geodesic arcs. Let $\underline{P}=\{P_i\}$ be a collection of marked points with $P_i\in\eta_i$. Let $\underline{w}=\{w_i\}$ be a collection of non-negative real numbers. The  \emph{strip deformation} of $Y$ determined by the triple $(\underline{\eta},\underline{P},\underline{w})$ is a new complete hyperbolic surface $Y'$ obtained by cutting along each $\eta_i$  and gluing in a strip of width $w_i$ so that the two points on the boundary of $Y\backslash\eta_i$  corresponding to the marked point $P_i\in\eta_i$ are identified with the feet of the strip. 

The length of every closed geodesic on $Y'$ is bigger than the corresponding length on $Y$, and strictly bigger if the geodesic in consideration intersects $\eta$ (see  \cite[Lemma 2.2]{DGK2016} for instance). Furthermore, if the weight $w_i\in \underline{w}$ goes to infinity, then the length of every closed geodesic intersecting the arc $\eta_i$ also goes to infinity. 

 Recall that every boarded hyperbolic surface is the convex hull of a unique complete hyperbolic surface. Whenever we say a strip deformation of a boarded hyperbolic surface $Y_0$, we mean the convex hull of the corresponding strip deformation of the complete hyperbolic surface which contains $Y_0$ as its convex hull. 

	\medskip
	
	In the remainder of the paper, we make the following convention. 
	\begin{quote}\label{con1}
		\begin{mdframed}
  \begin{convention}
  	We will also call a homologically trivial simple closed geodesic a \emph{separating} simple closed geodesic, and call a homologically nontrivial simple closed geodesic a \emph{nonseparating} simple closed geodesic. 
  \end{convention}
   \end{mdframed}
	\end{quote}
\begin{lemma}\label{lem:decreasing}
	Let $X\in\T_g$. Let $\{\gamma_1,\cdots,\gamma_k\}$ be a collection of disjoint simple closed geodesics. Let $\delta$ be a constant such that $\delta\geq \ell_X(\gamma_i)$. Then there exists a hyperbolic surface $X'\in\T_g$ such that 
	\begin{itemize}
		\item $\ell_{X'}(\gamma_i)=\delta$ for each $i$, and 
		\item $\ell_{X'}(\alpha)\geq \ell_X(\alpha)$ for every simple closed geodesic $\alpha$ on $X\backslash(\gamma_1\cup\cdots\cup\gamma_k)$. 
	\end{itemize}
\end{lemma}
\begin{proof}
	Let $X_1,\cdots,X_{m}$ be the components of $X\backslash(\gamma_1\cup\cdots\gamma_k)$. Applying strip deformations  to each $X_i$, we see that there exists a boarded hyperbolic surface $X_i'$ homeomorphic to $X_i$ such that
		\begin{itemize}
		\item the length of each geodesic boundary component of $X_i$ is exactly $\delta$,
		\item the length of every simple closed geodesic $\alpha$ in the interior of $X_i$ satisfes: $\ell_{X_i'}(\alpha)\geq \ell_{X_i}(\alpha)$.
	\end{itemize}
	Let $X'$ be a hyperbolic surface obtained by gluing $X_1',X_2',\cdots,X_{m}'$ along pairs of geodesic boundary components corresponding to $\gamma_i$ for all $i$. (The gluing process is not unique. The resulting surfaces differ by a twist about $\{\gamma_1,\cdots,\gamma_k\}$.)  Then $\ell_{X'}(\alpha)\geq \ell_X(\alpha)$ for every simple closed geodesic $\alpha$ on $X\backslash(\gamma_1\cup\cdots\cup\gamma_k)$ .  \end{proof}

\begin{remark}
There are other approaches to lengthen simple closed geodesics on boarded hyperbolic surfaces: using \emph{Fenchel-Nielsen coordinates} (see 
	 \cite{Parlier2005} for instance, as well as the following section), taking the \emph{Nielsen kernel} of a boarded surface \cite{Bers1976},  or more generally, using  \emph{pruning} constructions \cite{Dumas2006} (the inverse construction of grafting).  One advantage of strip deformation is that we could \emph{uniformly} lengthen all closed geodesics if the set $\underline{\alpha}$ of arcs cuts the underlying hyperbolic surface into topological disks.
\end{remark}

\section{Thick-thin decomposition and Fenchel-Nielsen coordinates}

The goal of this section is to use the Fenchel-Nielsen coordinates to lengthen short curves to some uniform constant while keeping other simple closed geodesics from being too contracted (Lemma \ref{lem:auxiliary:surface}). This will be used later for estimating the lower bound of the minimum of the algebraic intersection form in the moduli space. 
\subsection{Collar neighbourhoods and thick-thin decomposition of hyperbolic surfaces}\label{sec:thick:thin}
Given a hyperbolic surface $X\in\T_g$ and a simple closed geodesic $\gamma$, let $\C(\gamma)$ be the collar neighborhood of $\gamma$ defined by:		\begin{equation*}
			\C(\gamma):=\{p\in X:\mathrm{dist}(p,\gamma)\leq 
			\mathrm{arcsinh} \frac{1}{\sinh (\ell_X(\gamma_i)/2)}\}. 
		\end{equation*}
		Each circle inside $\C(\gamma)$ of distance $\rho$ from $\gamma$ has length 
		\begin{equation*}
			\ell_X(\gamma)\cosh \rho\leq \ell_X(\gamma) \sqrt{1+\frac{1}{\sinh^2 \ell_X(\gamma)/2}}\leq 2 \sqrt{1+\sinh^2 (\ell_X(\gamma)/2)}.
		\end{equation*}
		In particular, if  $\ell_X(\gamma)\leq 2\arcsinh1$, then each of the boundary components of $ \C(\gamma)$ has length at most $2\sqrt{2}$. Note that $\C(\gamma)$ is also foliated by geodesic arcs, each of which connects one boundary component to the other and intersects $\gamma$ perpendicularly. We call any such geodesic arc a \emph{perpendicular} of $\C(\gamma)$. Let $\xi$ be a geodesic arc from one boundary component of $\C(\gamma)$ to the other component. Suppose that $\xi$ crosses some perpendicular  $m$ times. Then it crosses any other perpendicular at most $m+1$ times. By an elementary calculation (\cite[Formula Glossary 2.2.2 (i)]{Buser2010}), we see that, if $\ell_X(\gamma)\leq 2\arcsinh1$, then
	\begin{equation}\label{eq:thin:length}
		|\ell_X(\xi)-(2|\log \ell_X(\gamma)|+m\ell_X(\gamma))|\leq C
	\end{equation} 
	for some universal constant $C$. 
	
		Hyperbolic surfaces with short curves admit  \emph{thick-thin decompositions}.
	Let $\{\gamma_1,\cdots,
	\gamma_k\}$ be the set of simple closed geodesics on $X$ of length at most $2\arcsinh1$. Then (i) $\gamma_1,\cdots,\gamma_k$ are pairwise disjoint, and (ii) any point outside $\cup_i \C(\gamma_i)$ has injective radius at least $\arcsinh~ 1\approx0.88$ \cite[Theorem 4.1.6]{Buser2010}. The union $\cup_i \C(\gamma_i)$ is called the \emph{thin} part of $X$ and its complementary is called the \emph{thick} part of $X$.  
	
	We end up this subsection with the following convention: 
	\begin{quote}
		\begin{mdframed}
		\begin{convention}\label{cvt:1}
			Fixed a constant $0<\consta\label{a:initial}<\arcsinh1\approx0.88$ once for all such that for any $0<x<\theconsta$, 
	\begin{equation}\label{cst:c}
		\arcsinh  \frac{1}{\sinh (x/2)}\geq |\log x|\geq \sqrt{2}.
	\end{equation}
		We may choose $\theconsta=e^{-\sqrt{2}}\approx0.24$.
		\end{convention}
	\end{mdframed}
	\end{quote}

\subsection{Fenchel-Nielsen coordinates}
\label{sec:fn}
In this subsection, we introduce the Fenchel-Nielsen coordinates on the Teichm\"uller space $\T_g$.  Let $\Gamma$  be a \textit{pants decomposition} of $S_g$ whose complementary $S\backslash \Gamma$ is a union of pairs of pants.

\begin{figure}
	\begin{tikzpicture}
		\draw(-3,2).. controls (-1,1) and (1,1).. (3,2); 
		\draw (-3,-2).. controls (-1,-1) and (1,-1).. (3,-2);
		
		\draw (-3.3,1).. controls (-2,0.5) and (-2,-0.5).. (-3.3,-1);   
		\draw (-3.3,1).. controls (-2.9,1.2) and (-2.8,1.6).. (-3,2);
		\draw (-3.3,1).. controls (-3.4,1.2) and (-3.4,1.6).. (-3,2);
		\draw (-3,1-3).. controls (-2.9,1.2-3) and (-2.8,1.6-3).. (-3.3,2-3);
		\draw (-3,1-3).. controls (-3.4,1.2-3) and (-3.4,1.6-3).. (-3.3,2-3); 
		
		\draw (3.3,1).. controls (2,0.5) and (2,-0.5).. (3.3,-1);   
		\draw (3.3,1).. controls (2.9,1.2) and (2.8,1.6).. (3,2);
		\draw (3.3,1).. controls (3.4,1.2) and (3.4,1.6).. (3,2);
		\draw (3,1-3).. controls (2.9,1.2-3) and (2.8,1.6-3).. (3.3,2-3);
		\draw (3,1-3).. controls (3.4,1.2-3) and (3.4,1.6-3).. (3.3,2-3); 
		
		\draw (0,1.25).. controls (-0.4, 1) and (-0.4,-1).. (0,-1.25);    
		\draw[dotted] (0,1.25).. controls (0.4, 1) and (0.4,-1).. (0,-1.25); 
		\draw (-0.4,-0.3)--(-0.3,-0.1)--(-0.2,-0.3);    
		
		\draw[red,line width=0.8pt](-2.93,-1.6)..controls (-2,-1.3) and (-1,-1).. (-0.21,-0.8);  
		\draw[blue,line width=0.8pt](2.93,1.6)..controls (2,1.3) and (1,1).. (-0.21,0.8);
		\draw[fill] (-0.22,0.8) circle (1pt) (-0.22,-0.8) circle (1pt); 
		\draw (-1.5,-0.9) node{$\eta'$} (-3.6,-1.5) node{$\gamma'$} (-0.5,0) node{$\gamma$}  (-1.6,0) node{$R'$}  (1.5,0) node{$R''$}  (1.5,0.9) node{$\eta''$} (3.6,1.5) node{$\gamma''$} ;
	\end{tikzpicture}
	\caption{\small Twists at $\gamma$ with respect to two pairs of pants $R'$ and $R''$.}
	\label{fig:FN}
\end{figure}
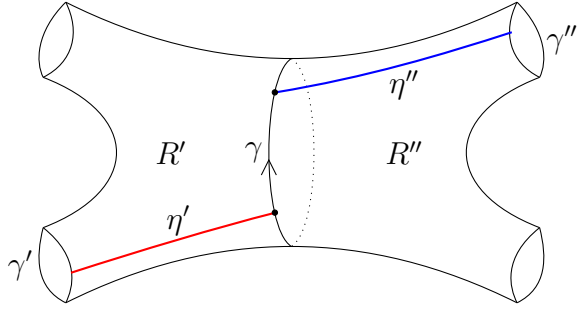

For $\gamma\in\Gamma$, the twist of $X\in\T_g$ with respect to $\Gamma$ at $\gamma$ is defined in the following way. Let $R'$ and $R''$ be the pairs of pants in $S\backslash\Gamma$, each of which contains $\gamma$ as a boundary component. Let $\gamma'$ (resp. $\gamma''$) be a geodesic boundary component of $R'$ (resp. $R''$) different from $\gamma$. Let $\eta'$ be the unique simple geodesic arc on $R'$ connecting $\gamma'$ and $\gamma$ perpendicularly. Let $\eta''$ be similarly defined. Choose an orientation of $\gamma$ such that $R'$ sits to the left of $\gamma$.  Let $d_{\gamma}(X)\in [-\ell_X(\gamma)/2,\ell_X(\gamma)/2)$ be the oriented distance from the foot $\eta'\cap\gamma$ to the foot $\eta''\cap\gamma$ along $\gamma$ with respect to prescribed orientation (see Figure \ref{fig:FN}). (In fact, the signed distance is independent of the choice of the orientation of $\gamma$. Instead, it depends on the orientation of the underlying surface.)   Let  
\begin{equation*}
	\widehat t_{\gamma}(X):=\frac{d_{
	\gamma}(X)}{\ell_X(\gamma)}\in [-\frac{1}{2},\frac{1}{2}).
\end{equation*}
The function $\widehat t_{\gamma}$ defines a function from $\T_g$ to the unit circle $\mathbb \R/\Z$.  Accordingly, the function $ \widehat t_{\Gamma}:=(\widehat t_{\gamma})_{\gamma\in\Gamma}$ defines a map from the Teichm\"uller space to the torus $\mathbb T^{3g-3}:=\R^{3g-3}/\Z^{3g-3}$. Let $X_0\in\T_g$ be a hyperbolic surface such that $\widehat t_{\Gamma}(X_0)=0$. Let $t_{\Gamma}:=\{t_{\gamma}\}_{\gamma\in\Gamma}: \T^{3g-3}\to \R^{3g-3}$ be the lift of $\widehat t_{\Gamma}:\T^{3g-3}\to \R^{3g-3}/\Z^{3g-3}$ such that $t_{\Gamma}(X_0)=0$.

The $(6g-6)$-tuple $(\ell_X(\gamma),t_\gamma(X))_{\gamma\in \Gamma}$ defines a Fenchel-Nielsen coordinates on the Teichm\"uller space $\T_g$.

  \begin{remark}\label{remark:twist} As we can see from the construction above, the Fenchel-Nielsen coordinate is not unique. The coordinate depends on the choice of $\Gamma$, the choice of $\eta'$ and $\eta''$ (which in turn is determined by the pairs $\{\gamma',\gamma\}$ and $\{\gamma'',\gamma\}$) for each $\gamma\in\Gamma$, and the choice of $X_0$. 
      \end{remark} 
      
      \subsection{Twists of proper arcs on pairs of pants}\label{sec:twist:arc}
      
      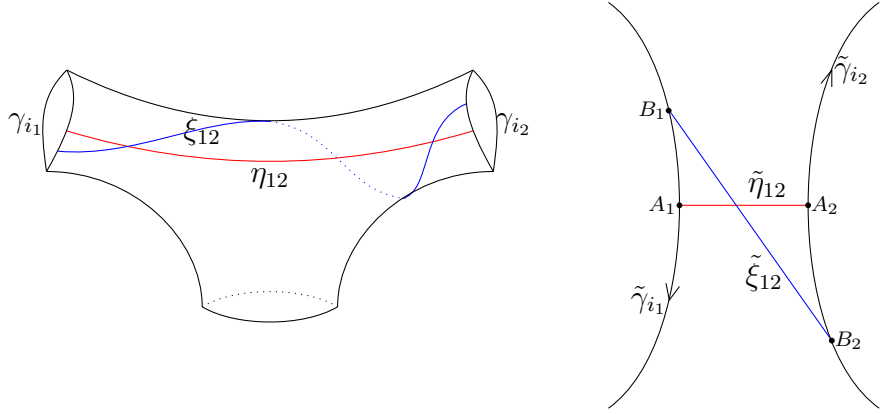
\begin{figure}
      	\begin{tikzpicture}[scale=0.9]
      		\draw(-3,2).. controls (-1,1) and (1,1).. (3,2); 
      			\draw[red](-3,1.1).. controls (-1,0.5) and (1,0.5).. (3,1.1);
      			\draw[blue](-3.13,0.8).. controls (-2,0.7) and (-1,1.3).. (0,1.24);
      			\draw[dotted, blue](0,1.24)..controls (1,1.1) and (1.1, 0.2).. (2,0.1);
      			\draw[blue](1.95,0.1).. controls (2.4, 0.2) and (2.2, 1.2).. (2.9,1.5);
		\draw[dotted] (-1,-1.5).. controls (-0.5,-1.2) and (0.5,-1.2).. (1,-1.5);
		\draw (-1,-1.5).. controls (-0.5,-1.8) and (0.5,-1.8).. (1,-1.5);
		
		\draw (-3.3,0.5).. controls (-2,0.5) and (-1,-0.5).. (-1,-1.5);   
		\draw (-3.3,0.5).. controls (-2.9,1.2) and (-2.8,1.6).. (-3,2);
		\draw (-3.3,0.5).. controls (-3.4,1.2) and (-3.4,1.6).. (-3,2);
		
		\draw (3.3,0.5).. controls (2,0.5) and (1,-0.5).. (1,-1.5);   
		\draw (3.3,0.5).. controls (2.9,1.2) and (2.8,1.6).. (3,2);
		\draw (3.3,0.5).. controls (3.4,1.2) and (3.4,1.6).. (3,2);

     \draw (0,0.4) node{$\eta_{12}$}
      (-1,1.1) node{$\xi_{12}$}
          (-3.6,1.2) node{$\gamma_{i_1}$}   (3.6,1.2) node{$\gamma_{i_2}$};
          
          \draw[red] (6.06,0)--(7.94,0);
          \draw[blue](5.9,1.4)--(8.3,-2);
         \draw (5.85,-1.15)--(5.9,-1.4)--(6.05,-1.2)  (8.14,1.8)--(8.3,2)--(8.3,1.75);
          \draw (5,3).. controls (6.4,2) and (6.4,-2).. (5,-3) (9,3).. controls (7.6,2) and (7.6,-2).. (9,-3);
          \draw (5.8,0) node {\tiny $A_1$}(8.2,0) node{\tiny $A_2$} (5.65,1.4)node{\tiny $B_1$} (5.6,-1.4)node{$\tilde\gamma_{i_1}$} (8.55,-2)node{\tiny $B_2$} (8.6,2)node{ $\tilde\gamma_{i_2}$} (7.3,-1)node{ $\tilde\xi_{12}$} (7.35,0.3)node{ $\tilde\eta_{12}$};
             \draw[fill] (6.05,0) circle(1pt) (7.94,0)
             (7.95,0) circle(1pt)  (5.89,1.4)circle (1pt) (8.3,-2)circle(1pt);
          
        		\end{tikzpicture}
		\caption{Twists of arcs in a pair of pants $P$.}
		\label{fig:twist}
      \end{figure}
      
      To estimate the length of simple closed geodesics in terms of Fenchel-Nielsen coordinates, we also need the \emph{twist} of proper simple arcs on pairs of pants with respect to the corresponding pants curves. More precisely, let $P$ be a pair of pants of $X\backslash\Gamma$ with boundary components, say  $\gamma_1,\gamma_2$, and $\gamma_3$. 
      
      Let $\xi_{12}$ be a simple arc on $P$ which connects $\gamma_{i_1}$ to $\gamma_{i_2}$.  Let $\eta_{12}$ be the unique geodesic homotopic to $\xi_{12}$ relative to $\gamma_{i_1}$ and $\gamma_{i_2}$, that is, the unique simple geodesic arc which connects $\gamma_{i_1}$ and $\gamma_{i_2}$ perpendicularly (Figure \ref{fig:twist} (a)).  Now, consider the lift of $P$ to the hyperbolic plane. Choose one lift $\tilde\eta_{12}$ of $\eta_{12}$ and let $\tilde{\gamma}_{i_1}$, $\tilde\gamma_{i_2}$, and $\tilde\xi_{12}$ be the corresponding lifts of $\gamma_{i_1}$, $\gamma_{i_2}$ and $\xi_{12}$. Let $A_j$ (resp. $B_j$) be the intersection point between $\tilde\gamma_{i_j}$ and $\tilde\eta_{12}$ (resp. $\tilde\xi_{12}$). Choose an orientation for  $\gamma_{i_j}$ such that $P$ sits to left  of  $\gamma_{i_j}$.  This induces an orientation on $\tilde\gamma_{i_j}$.  (see Figure \ref{fig:twist} (b)).  Let $\rho_{j}(\xi_{12})$  be the oriented distance from $A_j$ to $B_j$. The twist $t_j(\xi_{12})$ of $\xi_{12}$ at $\gamma_{i_j}$ is defined to be:
      \begin{equation}\label{eq:twist:xi}
      	t_j(\xi_{12}):=\frac{\rho_j(\xi_{12})}{\ell_X(\gamma_{i_j})}.
      \end{equation}
      By the hyperbolic trigonometry formula (\cite[Page 38, Equation (2.3.2)]{Buser2010}), we see that the length $\ell_X(\xi_{12})$ of the unique geodesic representative on $P$ which is homotopic to $\xi_{12}$ relative to the endpoints satisfies:
      \begin{eqnarray}
      &&\nonumber	\cosh \ell_X(\xi_{12})\\
      &=&\cosh \rho_1(\xi)\cosh \rho_2(\xi)\cosh \ell_X(\eta_{12})-\sinh\rho_1(\xi)\sinh\rho_2(\xi). \label{eq:dist:xi}
      \end{eqnarray}  
           About the length $\ell_X(\eta_{12})$ of $\eta_{12}$, if $\gamma_{i_1}\neq \gamma_{i_2}$, then, by \cite[Formula Glossary (page 454), 2.4.1(iii)]{Buser2010},  
      	\begin{equation}\label{eq:eta:12}
      		\cosh\ell_X(\eta_{12})=\frac{\cosh \frac{\ell_X(\gamma_{i_3})}{2}+\cosh \frac{\ell_X(\gamma_{i_1})}{2}\cosh \frac{\ell_X(\gamma_{i_2})}{2}}{\sinh\frac{\ell_X(\gamma_{i_1})}{2}\sinh \frac{\ell_X(\gamma_{i_2})}{2} },
      	\end{equation}
      	where $\gamma_{i_3}$ is the third  boundary component of $P$. 	
If $\gamma_{i_1}= \gamma_{i_2}$, let $\gamma_{i_3}$ be one of the other two boundary components of $P$ and let $\eta_{13}$ be the unique simple geodesic arc connecting $\gamma_{i_1}$ and $\gamma_{i_3}$ perpendicularly. By   \cite[Formula Glossary (page 454) 2.3.4 (i)]{Buser2010} we see that 
      	\begin{equation}\label{eq:eta:11}
      		\cosh \frac{\ell_X(\eta_{12})}{2}=\sinh \ell_X(\eta_{13})\sinh \frac{\ell_X(\gamma_{i_3})}{2}. 
      	\end{equation}

\subsection{Surfaces with short curves}\label{sec:short}
Given a hyperbolic surface $X\in\T_g$ with systolic length $\sys(X)< \ref{a:initial}$, we shall construct a new hyperbolic surface with big systole and controlled geometry. More precisely, let $\gamma_1,\cdots,\gamma_k$ be the set of simple closed geodesics on $X$ of length strictly less than $\ref{a:initial}$. In particular, they are pairwise disjoint \cite[Theorem 4.1.6]{Buser2010}.  We extend $\gamma_1,\cdots,\gamma_k$ to a pants decomposition $\Gamma:=\{\gamma_1,\cdots,\gamma_k,\gamma_{k+1}\cdots,\gamma_{3g-3}\}$. Let $(\ell_X(\gamma),t_\gamma(X))_{\gamma\in\Gamma}$ be the Fenchel-Nielsen coordinates defined in Section \ref{sec:fn}. 

Let $Y\in\T_g$ be the hyperbolic surface with Fenchel-Nielsen coordinates:
   \begin{equation*}
   \left(\epsilon_M,\cdots,\epsilon_M,\ell_X(\gamma_{k+1}),\cdots, \ell_X(\gamma_{3g-3}),t_{\gamma_1}(X),t_{\gamma_2}(X),\cdots,t_{\gamma_{3g-3}}(X) \right).
   \end{equation*}
In other words, we replace the length of each of $\gamma_1,\cdots,\gamma_k$ by $\ref{a:initial}$ while keeping the other terms unchanged.

\begin{lemma}\label{lem:auxiliary:surface}
	Let $X\in\T_g$ be a hyperbolic with $\sys(X)<\ref{a:initial}$. Let $Y$ be defined as above.  Then $\sys(Y)=\ref{a:initial}$.  Furthermore, for every simple closed geodesic $\alpha$ on $X$, the following holds.
	\begin{itemize}
		\item If $\alpha$ intersects the some $\gamma_i$,  then
			\begin{equation*}
	\frac{	\ell_{X}(\alpha)}{\ell_Y(\alpha)}\leq\ref{a:fn} \cdot  |\log \sys (X)|
	\end{equation*}
	for some universal constant $\consta\label{a:fn}$.
	 \item If $\alpha$ is disjoint from any of $\gamma_i$, then 
	 	\begin{equation*}
	\frac{	\ell_{X}(\alpha)}{\ell_Y(\alpha)}\leq 1.
	\end{equation*}
	\end{itemize}
	\end{lemma}
\begin{proof}

Let $\alpha$ be an arbitrary simple closed curve on $X$. Let $\alpha^Y$ be the geodesic on $Y$ which is homotopic to $\alpha$. The (geodesic representative of the) pants decomposition $\Gamma$ cuts $\alpha^Y$ into a collection $\{\xi_i^Y\}$ of geodesic segments.  From the construction of Fenchel-Nielson coordinates in Section \ref{sec:fn}, we see that $\alpha^Y$ (not the homotopy class of $\alpha^Y$) remains to be a simple closed curve on $X$, which is not necessarily a geodesic. To bound the length  $\ell_X(\alpha)$ of its geodesic representative on $X$ from above, we shall estimate the length $\mathrm{Length}_X(\xi^X_i)$ of the geodesic $\xi_i^X$ on $X$ which is homotopic to $\xi_i^Y$ relative to its endpoints:
\begin{equation*}
	\ell_X(\alpha)\leq \mathrm{Length}_X(\alpha^X)\leq \sum_i \mathrm{Length}_X(\xi_i^X).
\end{equation*}
To prove the lemma, it suffices to prove the following claim:

\begin{quote}
	\begin{mdframed}
	\textbf{Claim:} Let $\xi_i^X$ and $\xi_i^Y$ be as above.
		\begin{itemize}
		\item If $\xi_i^X$ is disjoint from any of the short curves $\gamma_1,\cdots,\gamma_k$, then \begin{equation*}
		\frac{\mathrm{Length}_X(\xi_i^X)}{\mathrm{Length}_Y(\xi_i^Y)}\leq 1.
	\end{equation*} 
	\item If $\xi_i^X$ intersects at least one of the short curves $\gamma_1,\cdots,\gamma_k$, then 	\begin{equation*}
		\frac{\mathrm{Length}_X(\xi_i^X)}{\mathrm{Length}_Y(\xi_i^Y)}\leq \ref{a:fn}\cdot |\log \sys(X)|
	\end{equation*} 
	for some universal constant $\ref{a:fn}$.
		\end{itemize}
	\end{mdframed}
\end{quote}
 Indeed, if $\alpha$ does not intersect any of the short curves $\gamma_1,\cdots,\gamma_k$, then by the first item of the claim, we see that
 \begin{equation*}
	\ell_X(\alpha)\leq \mathrm{Length}_X(\alpha^X)\leq \sum_i \mathrm{Length}_X(\xi_i^X)\leq  \sum_i \mathrm{Length}_Y(\xi_i^Y)=\ell_Y(\alpha).
\end{equation*}
If $\alpha$ intersects at least one of the short curves $\gamma_1,\cdots,\gamma_k$, then by the second item of the claim, we see that
 \begin{eqnarray*}
	\ell_X(\alpha)&\leq &\mathrm{Length}_X(\alpha^X)\leq \sum_i \mathrm{Length}_X(\xi_i^X)\\
	&\leq&  \ref{a:fn}|\log \sys(X)|\cdot  \sum_i \mathrm{Length}_Y(\xi_i^Y)=\ref{a:fn}|\log \sys(X)|\cdot \ell_Y(\alpha).
\end{eqnarray*}
It remains to consider the systolic length $\sys(Y)$. Let $\alpha$ be simple closed geodesic on $Y$. If $\alpha$ is one of $\gamma_i$, then $\ell_X(\gamma_i)=\ref{a:initial}$. If $\alpha$ is disjoint from any of $\gamma_i$, then the first item in the claim implies that $\ell_Y(\alpha)\geq \ell_X(\alpha)\geq \ref{a:initial}$. If $\alpha$ intersects some $\gamma_i$. then it crosses the collar neighborhood $\C(\gamma_i)$ of $\gamma_i$, which is of width 
$\arcsinh(1/\sinh(\ell_X(\gamma)/2))\geq \ref{a:initial}$ due to the assumption that $\ell_X(\gamma_i)<\ref{a:initial}$. Hence, $\ell_Y(\alpha)\geq 2\ref{a:initial}$. This proves that $\sys(Y)\geq \ref{a:initial}$. That $\ell_X(\gamma_i)=\ref{a:initial}$ for $i\leq k$ implies that $\sys(Y)=\ref{a:initial}$.

The remainder of the proof is devoted to proving the claim. To simplify notations, we set $\ell_X(\xi^X_i):=\mathrm{Length}_X(\xi^X_i)$ and $\ell_Y(\xi^Y_i):=\mathrm{Length}_Y(\xi^Y_i)$.  Let $P$ be the pair of pants which contains $\xi_i^X$. Let $\gamma_{i_1}$ and $\gamma_{i_2}$ be two pants curves connected by the arc $\xi^X_i$ (and hence also by $\xi^Y_i$). Let $t_{1}(\xi)$ (resp. $t_{2}(\xi)$) be the twist of $\xi^X$ at $\gamma_{i_1}$ (resp. $\gamma_{i_2}$) with respect to  the pair of pants $P$ (see Section \ref{sec:twist:arc} for relevant definitions). By the definition of $\xi^Y$, we see that $t_{1}(\xi)$ (resp. $t_{2}(\xi)$) is also the twist of $\xi^Y$ at $\gamma_{i_1}$ (resp. $\gamma_{i_2}$) with respect to  the pair of pants $P$.
Let $\eta_{i}^X$ be the geodesic segment under the metric $X$  which is homotopic to $\xi_i^X$  relative to $\gamma_{i_1}$ and $\gamma_{i_2}$.   Let $\eta_{i}^Y$ be the geodesic segment under the metric $Y$  which is homotopic to $\xi_i^Y$  relative to $\gamma_{i_1}$ and $\gamma_{i_2}$. There are two case depending on whether $\gamma_{i_1}=\gamma_{i_2}$ or not.

\textbf{Case 1. $\gamma_{i_1}\neq\gamma_{i_2}$.} Let $\gamma_{i_3}$ be the third edge of the pair of pants $P$ containing $\xi^X$. There are two more subcases depending on whether $\gamma_{i_1}$ or $\gamma_{i_2}$ is short.

\underline{Subcase 1-1:  neither $\gamma_{i_1}$ nor $\gamma_{i_2}$ is short}, that is $\ell_X(\gamma_{i_1})\geq \ref{a:initial}$ and $\ell_X(\gamma_{i_2})\geq \ref{a:initial}$. Then, by the construction of $Y$,
\begin{equation}\label{eq:gamma:equal}
	\ell_X(\gamma_{i_1})= \ell_Y(\gamma_{i_1}) \text{ and }\ell_X(\gamma_{i_2})= \ell_Y(\gamma_{i_2}).
\end{equation}
Since $\ell_X(\gamma_{i3})\leq \ell_Y(\gamma_{i3})$, it follows from \eqref{eq:eta:12} that $\ell_X(\eta_i^X)\leq \ell_Y(\eta_i^Y)$. Hence, 
by \eqref{eq:dist:xi}, we see that 
$\ell_X(\xi^X_i)\leq \ell_Y(\xi^Y_i)$. 

\underline{Subcases 1-2: either $\gamma_{i_1}$ or  $\gamma_{i_2}$ is short,} say $\ell_X(\gamma_{i_1})<\ref{a:initial}$. Accordingly, we have $\ell_Y(\gamma_{i_1})=\ref{a:initial}$. By \eqref{eq:eta:12}, we see that 
\begin{eqnarray}
	\nonumber\ell_X(\eta_i^X)&\geq& \log \frac{\cosh \frac{\ell_X(\gamma_{i_1})}{2}\cosh \frac{\ell_X(\gamma_{i_2})}{2}}{\sinh \frac{\ell_X(\gamma_{i_1})}{2}\sinh \frac{\ell_X(\gamma_{i_2})}{2}}\\
	&\geq&\log \frac{\cosh \frac{\ell_X(\gamma_{i_1})}{2}}{\sinh \frac{\ell_X(\gamma_{i_1})}{2}}\geq 
	\log \frac{\cosh \frac{\ref{a:initial}}{2}}{\sinh \frac{\ref{a:initial}}{2}}, \label{eq:eta:X:12}
	\end{eqnarray}
and similarly, 
\begin{equation}\label{eq:eta:Y:12}
	\ell_Y(\eta_i^Y)\geq \log \frac{\cosh \frac{\ref{a:initial}}{2}}{\sinh \frac{\ref{a:initial}}{2}}.
	\end{equation}
Again by \eqref{eq:eta:12}, we have 
\begin{eqnarray*}
	&&\frac{\cosh\ell_X(\eta_i^X)}{\cosh \ell_Y(\eta_i^Y)}\\&=&\frac{\cosh \frac{\ell_X(\gamma_{i_3})}{2}+\cosh \frac{\ell_X(\gamma_{i_1})}{2}\cosh \frac{\ell_X(\gamma_{i_2})}{2}}{\cosh \frac{\ell_Y(\gamma_{i_3})}{2}+\cosh \frac{\ell_Y(\gamma_{i_1})}{2}\cosh \frac{\ell_Y(\gamma_{i_2})}{2}}\cdot \frac{\sinh\frac{\ell_Y(\gamma_{i_1})}{2}\sinh \frac{\ell_Y(\gamma_{i_2})}{2} }{\sinh\frac{\ell_X(\gamma_{i_1})}{2}\sinh \frac{\ell_X(\gamma_{i_2})}{2} }\\ 
	&\leq & \frac{\sinh\frac{\ell_Y(\gamma_{i_1})}{2}\sinh \frac{\ell_Y(\gamma_{i_2})}{2} }{\sinh\frac{\ell_X(\gamma_{i_1})}{2}\sinh \frac{\ell_X(\gamma_{i_2})}{2} }\leq  \frac{\sinh ^2\frac{\ref{a:initial}}{2}}{\sinh^2\frac{\sys(X)}{2}}
\end{eqnarray*}
where in the last equality we use the assumption that $\ell_X(\gamma_{i_1})<\ref{a:initial}=\ell_Y(\gamma_{i_1})$ and the fact  that $\ell_Y(\gamma_{i_2})=\ell_X(\gamma_{i_2})$ if $\ell_X(\gamma_{i_2})\geq \ref{a:initial}$ or $\ell_Y(\gamma_{i_2})=\ref{a:initial} $ if  $\ell_X(\gamma_{i_2})<\ref{a:initial}$.
 Hence
\begin{eqnarray*}
	\ell_X(\eta_i^X)-\ell_Y(\eta_i^Y)&\leq& \log2+\log \frac{\cosh \ell_X(\eta_i^X)}{\cosh \ell_Y(\eta_i^Y)}\leq \log 2+2\log \frac{\sinh \frac{\ref{a:initial}}{2}}{\sinh\frac{\sys(X)}{2}}.
\end{eqnarray*}
Together with \eqref{eq:eta:Y:12}, this implies that
\begin{eqnarray}
\nonumber
	\frac{\ell_X(\eta_i^X)}{\ell_Y(\eta_i^Y)}&=& 1+\frac{\ell_X(\eta_i^X)-\ell_Y(\eta_i^Y)}{\ell_Y(\eta_i^Y)}
\\
	&\leq &  \ref{c:ratio:eta:12}
	\cdot |\log \sys(X)|\label{eq:ratio:eta:12} 
\end{eqnarray}
for some universal positive constant $\constc\label{c:ratio:eta:12}$.
Combining  \eqref{eq:twist:xi}, \eqref{eq:dist:xi}, \eqref{eq:eta:X:12}, and \eqref{eq:eta:Y:12}, we see that there exist two universal constants $\constc \label{c:xi:x}$ 
 and  $\constc \label{c:xi:y}$ such that \begin{eqnarray*}
	 \ell_X(\xi_i^X)&\leq& \ref{c:xi:x} \cdot ( t_{i_1}(\xi_i)\ell_X(\gamma_{i_1})+ t_{i_2}(\xi)\ell_X(\gamma_{i_2})+ \ell_X(\eta_i^X))
\end{eqnarray*}
and
\begin{eqnarray*}
	 	 \ell_Y(\xi_i^Y)
	 &\geq & \ref{c:xi:y}\cdot (t_{i_1}(\xi_i)\ell_Y(\gamma_{i_1})+t_{i_2}(\xi)\ell_Y(\gamma_{i_2})+\ell_Y(\eta_i^Y )).
\end{eqnarray*}
Therefore,
\begin{eqnarray*}
	\frac{\ell_X(\xi^X_i)}{\ell_Y(\xi^Y_i)}&\leq &\frac{\ref{c:xi:x}}{\ref{c:xi:y}} \cdot \frac{t_{i_1}(\xi_i)\ell_X(\gamma_{i_1})+ t_{i_2}(\xi)\ell_X(\gamma_{i_2})+ \ell_X(\eta_i^X)}{t_{i_1}(\xi_i)\ell_Y(\gamma_{i_1})+ t_{i_2}(\xi)\ell_Y(\gamma_{i_2})+ \ell_Y(\eta_i^Y)}\\
	&\leq & \frac{\ref{c:ratio:eta:12}\ref{c:xi:x}}{\ref{c:xi:y}}\cdot |\log \sys(X)|. \qquad\qquad (\text{ by }\eqref{eq:ratio:eta:12}) 
\end{eqnarray*}

\bigskip
\textbf{Case 2. $\gamma_{i_1}=\gamma_{i_2}$.}  Let $\gamma_{i_3}$ be one of the other two geodesic boundary components of the pair of pants $P$. Let $\eta_{i_{13}}^X$ (resp. $\eta_{i_{13}}^Y)$  be the unique simple geodesic arc connecting $\gamma_{i_1}$ and $\gamma_{i_3}$ perpendicularly on $X$ (resp. $Y$). Using \eqref{eq:eta:11} and \eqref{eq:eta:12}, we see that
\begin{eqnarray}
 &&\cosh^2\ell_X(\eta_i^X) \label{eq:cosh:eta:11:X}
	\\ \nonumber &=&\sinh^2 \frac{\ell_X(\gamma_{i_3})}{2}\sinh^2 \ell_X(\eta_{i_{13}}^X)\\
	\nonumber
	&=& \sinh^2 \frac{\ell_X(\gamma_{i_3})}{2} \left[\frac{\left(\cosh \frac{\ell_X(\gamma_{i_3'})}{2}+\cosh \frac{\ell_X(\gamma_{i_1})}{2}\cosh \frac{\ell_X(\gamma_{i_3})}{2}\right)^2}{\sinh^2\frac{\ell_X(\gamma_{i_1})}{2}\sinh ^2\frac{\ell_X(\gamma_{i_3})}{2} }-1\right]\\
	&=& \nonumber \sinh^{-2}\frac {\ell_X(\gamma_{i_1})}{2} \left(\cosh^2 \frac{\ell_X(\gamma_{i_3'})}{2}+\cosh^2 \frac{\ell_X(\gamma_{i_1})}{2}+\cosh^2 \frac{\ell_X(\gamma_{i_3})}{2}\right. \\
	&&\nonumber \qquad \left. + 2\cosh\frac{\ell_X(\gamma_{i_3'})}{2}\cosh \frac{\ell_X(\gamma_{i_1})}{2}\cosh \frac{\ell_X(\gamma_{i_3})}{2}-1\right).
\end{eqnarray}
 Similarly, we have
\begin{eqnarray}\label{eq:cosh:eta:11:Y}
&&	\cosh^2\ell_Y(\eta_i^Y)  \\
&=&  \nonumber \sinh^{-2}\frac {\ell_Y(\gamma_{i_1})}{2} \left(\cosh^2 \frac{\ell_Y(\gamma_{i_3'})}{2}+\cosh^2 \frac{\ell_Y(\gamma_{i_1})}{2}+\cosh^2 \frac{\ell_Y(\gamma_{i_3})}{2}\right. \\
	&&\nonumber \qquad \left. + 2\cosh\frac{\ell_Y(\gamma_{i_3'})}{2}\cosh \frac{\ell_Y(\gamma_{i_1})}{2}\cosh \frac{\ell_Y(\gamma_{i_3})}{2}-1\right).
\end{eqnarray}
Since $\ell_X(\gamma_{i})\leq \ell_Y(\gamma_i)$,  this implies that
\begin{eqnarray}\label{eq:ratio:cosh:square}
	\frac{\cosh\ell_X(\eta_i^X)}{\cosh\ell_Y(\eta_i^Y)}\leq \frac{\sinh\frac{\ell_Y(\gamma_{i_1})}{2}}{\sinh\frac{\ell_X(\gamma_{i_1})}{2}}.
\end{eqnarray}

\underline{Subcase 2-1: $\ell_X(\gamma_{i_1})\geq \ref{a:initial}$.}  In particular, $\ell_X(\gamma_{i_1})=\ell_Y(\gamma_{i_1})$. Hence by \eqref{eq:ratio:cosh:square}, we see that $\ell_X(\eta_i^X)\leq \ell_Y(\eta_i^Y)$. Combining with  \eqref{eq:dist:xi}, and the assumption that $\gamma_{i_1}=\gamma_{i_2}$, we infer that $\ell_X(\xi_i^X)\leq \ell_Y(\xi_i^Y)$.

\underline{Subcase 2-2: $\ell_X(\gamma_{i_1})< \ref{a:initial}$.}  Then $\ell_Y(\gamma_{i_1})= \ref{a:initial}$.  It follows from \eqref{eq:cosh:eta:11:X} and \eqref{eq:cosh:eta:11:Y} that
\begin{equation*}
	\ell_X(\eta_i^X)\geq \log \frac{2}{\sinh \frac{\ref{a:initial}}{2}}\quad\text{ and } \quad\ell_Y(\eta_i^Y)\geq \log \frac{2}{\sinh \frac{\ref{a:initial}}{2}}.
\end{equation*}
Furthermore,  from  \eqref{eq:ratio:cosh:square}, we see that
\begin{equation*}
	\frac{\cosh\ell_X(\eta_i^X)}{\cosh\ell_Y(\eta_i^Y)}\leq \frac{\sinh\frac{\ref{a:initial}}{2}}{\sinh\frac{\sys(X)}{2}}.
\end{equation*}
Applying an argument similarly as in the Subcase 1-2 yields 
\begin{equation*}
	\frac{\ell_X(\xi^X_i)}{\ell_Y(\xi^Y_i)}\leq \ref{c:ratio:eta:11}\cdot |\log \sys(X)|
\end{equation*}
for some constant universal constant $\constc\label{c:ratio:eta:11}$. 

In summary, the first item of the claim follows from the discussion in Subcase 1-1 and Subcase 2-1. The second item of the claim follows from the discussion in Subcase 1-2 and Subcase 2-2 by setting $\ref{a:fn}=\max\{\frac{\ref{c:ratio:eta:12}\ref{c:xi:x}}{\ref{c:xi:y}},\ref{c:ratio:eta:11}\}$. The proof of the claim, and hence of the lemma, is now complete. 
\end{proof}
\begin{remark}
	Instead of using strip deformations, we could also use the second item of Lemma \ref{lem:auxiliary:surface} or \cite{Parlier2005} to prove Lemma \ref{lem:decreasing}. The reason for using strip deformation is that strip deformation is intuitive and straightforward.
	\end{remark}

\section{Existence of minimum in moduli space}
The goal of this section is to prove Theorem \ref{thm:minimum}, which states that the algebraic intersection form $K(X)$ has a minimum in the moduli space.

We start the proof with the following lemma.

\begin{lemma}\label{lem:intersection:direction}
			Let $X\in\T_g$ be a hyperbolic surface with a simple closed geodesics $\gamma$ of length at most  $\ref{a:initial}$. If $\alpha$ is an  oriented simple closed geodesic with $|\Int(\gamma,\alpha)|<\i(\gamma,\alpha)$, then there exist oriented simple  closed geodesics $\alpha_1,\cdots,\alpha_k$ such that 
			\begin{enumerate}[(i)]
				\item  $\alpha_1\cup\alpha_2\cup\cdots\cup\alpha_k$ is homologous to  $\alpha$, and
				\item $\ell_X(\alpha_1)+\ell_X(\alpha_2)+\cdots+\ell_X(\alpha_k)\leq \ell_X(\alpha) $, and
				\item $|\Int(\gamma,\alpha_i)|=\i(\gamma,\alpha_i)$ for every $\alpha_i$.
			\end{enumerate}
		\end{lemma}
		\begin{proof}
		  Notice that $\gamma$ cuts $ \alpha$ into  $\i(\gamma,\alpha)$ segments.  Let $\alpha_{1}$ be a component of $\alpha\backslash\gamma$ such that the algebraic intersection from $\gamma$ (for any choice of orientation) to $\alpha$  has different signs at the endpoints of $\alpha_{1}$.  Notice that $\alpha_{1}$ cuts $\gamma$ into two subsegments. Choose one of them with appropriate orientation, say $\gamma_{1}$, such that $\alpha_{1}  \cup \gamma_{1}$ is an oriented simple closed curve.  Let $\widehat{\alpha}_{1}$ be a  simple closed geodesic homotopic to $\alpha_{1}\cup \gamma_{1}$.  Let $\gamma_{1}^-$ be the segment underlying $\gamma_{1}$ with \emph{opposite} orientation. Let $\widehat{\alpha}_{1}'$ be the oriented closed geodesic homotopic to $(\alpha\backslash\alpha_{1}) \cup \gamma_{1}^-$. From the construction, we see that 
	\begin{itemize}
		\item  $\widehat{\alpha}_{1}\cup \widehat{\alpha}_{1}'$  is homologous to $\alpha$, 
		\item $\i(\widehat{\alpha}_{1}\cup \widehat{\alpha}_{1}',\gamma)\leq\i(\alpha,\gamma)-2$.
	\end{itemize}
	   Moreover, the assumption that $\ell(\gamma)\leq  \ref{a:initial}$ implies that 
	\begin{equation}\label{eq:hat:alpha:short}
 \ell_X(\widehat{\alpha}_{1}) +\ell_X(\widehat{\alpha}_{1}') \leq \ell_X(\alpha).
	\end{equation}
	To see this, 	let $\xi$ be  the component $(\alpha_1\cup\gamma_1)\cap \C(\gamma)$ which contains $\gamma_1$. Then $\xi$ is homotopic to a subsegment $\xi_0$ of the boundary component of $\C(\gamma)$. The assumption $\ell(\gamma)\leq  \ref{a:initial}$ implies that $\gamma$ admits a collar neighbourhood $\C(\gamma)$ of width 
	\begin{equation*}
		\arcsinh\frac{1}{\sinh (\ell_X(\gamma)/2)}\geq \sqrt{2}.
	\end{equation*}
	Furthermore, each of the boundary components of $\C(\gamma)$ has length at most $2\sqrt{2}$ (see Section \ref{sec:thick:thin}).  
	Then the length $\mathrm{Length}_X(\xi\backslash\gamma_1)$ of $\xi\backslash\gamma_1$, which consists of two geodesic arcs connecting the boundary of $\C(\gamma)$ to $\gamma$,  and the length $\mathrm{Length}_X(\xi_0)$ of $\xi_0$ satisfy:
	\begin{equation*}
		\mathrm{Length}_X(\xi\backslash\gamma_1)\geq 2\sqrt{2}\geq \mathrm{Length}_X(\xi_0).
	\end{equation*}
	Hence, we could homotopy $\alpha_1\cup\gamma_1$ to a new curve by replacing $\xi$ with $\xi_0$  so that the resulting curve, and hence its geodesic representative $\widehat\alpha_1$, is of length less than $\ell_X(\alpha_1)$. Similarly, we have $\ell_X(\widehat\alpha_1')\leq \ell_X(\alpha\backslash\alpha_1)$, completing the proof of \eqref{eq:hat:alpha:short}.

	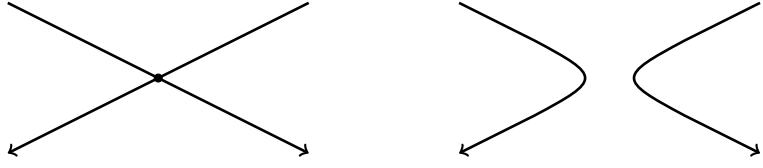
\begin{figure}
	 
	 \begin{tikzpicture}
	 	\draw[line width=1pt][->] (-2,1)--(2,-1);
	 	\draw[line width=1pt][->] (2,1)--(-2,-1);
	 	\draw[fill] (0,0) circle (1.5pt); 
	 	
	 	\draw[line width=1pt][->] (-2+6,1)--(-1+6,0.5)..controls  (-0.1+6, 0.02) and (-0.1+6,-0.02)..(-1+6,-0.5)--(-2+6,-1);
	 		\draw[line width=1pt][->] (2+6,1)--(1+6,0.5)..controls  (0.1+6, 0.02) and (0.1+6,-0.02)..(1+6,-0.5)--(2+6,-1);
	 \end{tikzpicture}
	 \caption{Smooth operation at a self-intersection.}
	 \label{fig:smooth}
	\end{figure}
	
	Note that $\widehat{\alpha}_{1}$ is  simple while $\widehat{\alpha}_1'$ is possibly not. If $\widehat{\alpha}_{1}'$ is not simple, then we apply a smoothing operation at each self-intersection point of $\widehat{\alpha}_1'$ (see Figure \ref{fig:smooth}). This gives us several oriented simple closed geodesics $\alpha_{11}',\cdots,\alpha_{1k}'$ such that
	\begin{itemize}
		\item $\widehat\alpha_{11}'\cup\cdots\cup\widehat\alpha_{1k}'$  is homologous to $\widehat{\alpha}_1' $,
		\item $\sum_j \i(\gamma,\widehat{\alpha}_{1j}')\leq \i(\gamma, \widehat\alpha_1')$,
		\item $\sum_j \ell_X(\widehat\alpha_{1j}')\leq \ell_X(\widehat\alpha_1')$.
	\end{itemize}
	
	Set $A_1:=\{\widehat{\alpha}_{1},\widehat{\alpha}_{11}',\cdots,\widehat{\alpha}_{1k}'\}$. Then, 
		\begin{itemize}
		\item  $\cup_{\alpha'\in A_1}{\alpha'}$ is homologous to $\alpha$,
		\item $\sum_{\alpha'\in A_1} \i(\alpha',\cup\gamma_i)\leq \i(\alpha,\gamma)-2$,
		\item $ \sum_{\alpha'\in A_1}\ell_X(\alpha')\leq \ell_X(\alpha)$.
	\end{itemize}
	 Inductively, suppose that $A_{i-1}$ has been constructed. We apply the above construction for each of $\alpha'\in A_{i-1}$ with $|\Int(\alpha',\gamma)|<\i(\alpha',\gamma)$. This yields a set $A_{i}$ of oriented simple closed geodesics such that 
	\begin{itemize}
		\item  $\cup_{\alpha'\in A_i}{\alpha'}$ is homologous to $\alpha$,
		\item $\sum_{\alpha'\in A_i} \i(\alpha',\cup\gamma_i)\leq \i(\alpha,\gamma)-2i$,
		\item $ \sum_{\alpha'\in A_i}\ell_X(\alpha')\leq \ell_X(\alpha)$.
	\end{itemize}
	Since for each $i$, $\sum_{\alpha'\in A_i} \i(\alpha',\cup\gamma_i)\geq |\Int(\alpha,\gamma)|$,  we see that the construction stops after $m\leq \frac{\i(\alpha,\gamma)-|\Int(\alpha,\gamma)|}{2}$ steps. The resulting set $A_m$ satisfies the required properties. 
		\end{proof}

Note that the defining equation  \eqref{eq:intection:length} of the algebraic intersection form $K(\cdot )$ is also well-defined for boarded surfaces. This allows us to extend the definition of $K(\cdot)$ to hyperbolic surfaces with geodesic boundary components. 

\begin{lemma}\label{lem:near:bnd}
	Let $X\in T_g$ be a hyperbolic surface with a collection of disjoint separating simple closed geodesics  $\gamma_1,\cdots,\gamma_k$ of length at most $\ref{a:initial}$.  Then,	\begin{equation*}
		K(X)=\max_{i} K(X_i)
	\end{equation*}
	where $X_1, \cdots, X_{k+1}$ are  components of $X\backslash(\gamma_1\cup\cdots\cup\gamma_k)$.
\end{lemma}
\begin{proof}
	Let $X_1,X_2,\cdots, X_{k+1}$ be the components of $X\backslash(\gamma_1\cup \cdots\cup \gamma_k)$.  Let $\alpha$ and $\beta$ be two simple closed geodesics. Since $\gamma_i$ is a separating simple closed geodesic, it follows that $\Int(\alpha,\gamma_i)=\Int(\beta,\gamma_i)=0$. Applying  Lemma \ref{lem:intersection:direction} to each $\gamma_i$ and $\alpha$,	we see that there exist $2k+2$  sets $\Omega_1,\cdots,\Omega_{k+1}$, $\Psi_1,\Psi_1,\cdots,\Psi_{k+1}$ of oriented simple closed geodesics  such that 
		\begin{itemize}
		\item $\cup_{\alpha'\in \Omega_i}\alpha'\subset X_i$ and $\cup_{\beta'\in \Psi}\beta'\subset X_i$ for each $i$
		\item  $\cup_{i}\cup_{\alpha'\in\Omega_i}\alpha'$   is homologous to $\alpha$ and $\cup_{i}\cup_{\beta'\in\Psi_i}\beta'$   is homologous to  $\beta$,
		\item $\sum_i\sum_{\alpha'\in\Omega_i}\ell_X(\alpha')\leq \ell_X(\alpha)$ and $\sum_i\sum_{\beta'\in\Psi_i}\ell_X(\beta')\leq \ell_X(\beta)$.
	\end{itemize}
	Then,
	\begin{eqnarray*}
		\frac{|\Int(\alpha,\beta)|}{\ell_X(\alpha)\ell_X(\beta)}&\leq&\frac{|\sum_i\Int(\cup_{\alpha'\in\Omega_{i}}\alpha',\cup_{\beta'\in\Psi_{i}}\beta')|}{(\sum_i\sum_{\alpha'\in\Omega_{i}}\ell_X(\alpha'))(\sum_i\sum_{\beta'\in\Psi{i}}\ell_X(\beta'))}\\
		&\leq &\max_i\frac{|\Int(\cup_{\alpha'\in\Omega_{i}}\alpha',\cup_{\beta'\in\Psi_{i}}\beta')|}{(\sum_{\alpha'\in\Omega_{i}}\ell_X(\alpha'))(\sum_{\beta'\in\Psi{i}}\ell_X(\beta'))}
		\\&\leq & \max_i\max\limits_{\alpha'\in\Omega_i,\beta'\in\Psi_i}\frac{|\Int(\alpha',\beta')|}{\ell_X(\alpha')\ell_X(\beta')}
		\\&\leq &
		  \max_i K(X_i).
	\end{eqnarray*}	
		Consequently, 
	\begin{equation*}
		K(X)=\sup_{\alpha,\beta}\frac{|\Int(\alpha,\beta)|}{\ell_X(\alpha)\ell_X(\beta)}\leq \max_i K(X_i).
	\end{equation*}
	\end{proof}

 We also need the following estimate, which is essentially the lower bound of \cite[Theorem 4.2]{MM2014}. For reader's convenience, we provide a proof.

\begin{lemma}\label{lem:nonsys}
Let $X\in\T_g$ be a hyperbolic surface with $\sys_h(X)\leq \ref{a:initial}$. Then,
	\begin{equation*}
		K(X)\geq \frac{1}{\ref{a:short:systole}\cdot  g\cdot \sys_h(X)(g+|\log \sys_h(X)|)}
	\end{equation*}
	for some universal positive constant $\consta\label{a:short:systole}$. 
\end{lemma}
\begin{proof}
	Let $\gamma$ be one of the shortest nonseparating simple closed geodesics. In particular, $\ell_X(\gamma)=\sys_h(X)$. By \cite[Theorem 5.1.4]{Buser2010}, there exists a pants decomposition $\Gamma$ which contains $\gamma$ such that
	\begin{equation}
		\ell_X(\gamma')\leq 6\sqrt{3\pi}(g-1)
	\end{equation}
	for every $\gamma'\in\Gamma$. Let $\beta$ be a simple closed geodesic that intersects $\gamma$ exactly once and intersects every other simple closed curve in $\Gamma$ at most once. Next, we shall estimate the length of $\beta$ from above.
	Let $P$ in $X\backslash\Gamma$ be an arbitrary pair of pants with boundary components $\gamma_i,\gamma_j,\gamma_k$. Let $\delta$ be  a geodesic arc  in $P$ perpendicular to $\gamma_i$ and $\gamma_j$ with $\gamma_i\neq\gamma_j$. Then by \eqref{eq:eta:12}, we see that 
	\begin{eqnarray*}
		\cosh \ell_X(\delta)&=&\frac{\cosh(\ell_X(\gamma_k)/2)+\cosh(\ell_X(\gamma_i)/2)\cosh(\ell_{X}(\gamma_j)/2}{\sinh(\ell_X(\gamma_i)/2)\sinh(\ell_{X}(\gamma_j)/2}\\
		&\leq& \frac{2\cosh^2(3\sqrt{3\pi}(g-1))}{\sinh^2(\sys_h(X)/2)}.
	\end{eqnarray*}
	Since $\sys_h(X)\leq \ref{a:initial}$, it follows that there exists a positive constant $\constc\label{c:pant:seam}$ such that $\sinh(\ell_X(\alpha)/2)\geq \ref{c:pant:seam}\ell_X(\alpha)$.
	Hence
	\begin{equation*}
		\ell_X(\delta)\leq \log2+\log\cosh \ell_X(\delta)\leq 4\log 2+6\sqrt{3\pi}(g-1)+2\log\frac{1}{\ref{c:pant:seam}}+2\log \frac{1}{\sys_h(X)}. 
	\end{equation*}
	Since $\beta$ intersects each of $\{\gamma_i\}$ at most once, we may homotopic $\beta$ to a new curve 
	\begin{equation*}
		\delta_1*\xi_1*\delta_2*\xi_2*\cdots *\delta_m*\xi_m
	\end{equation*}	
	with $m\leq |\Gamma|=3g-3$,
	where $\delta_i$ is a geodesic arc perpendicularly connecting distinct boundary components of some pair of pants in $X\backslash\Gamma$ and $\xi_i$ is a subarc of some of $\{\gamma_i\}$.  Combining with the estimate above about   the length of $\delta_i$, we see that
	\begin{equation}\label{eq:nonsys:beta}
		\ell_X(\beta)\leq \sum_i (\ell_X(\delta_i)+\ell_X(\xi_i))\leq \ref{c:beta:short}\cdot  g (g+\log\frac{1}{\sys_h(X)})
	\end{equation} 
	for some universal constant $\constc\label{c:beta:short}$.
	Combined with the definition of the algebraic intersection form, this implies that
	\begin{eqnarray*}
		K(X)&\geq& \frac{|\Int(\gamma,\beta)|}{\ell_X(\gamma)\ell_X(\beta)}
		\\
		&\geq &  \frac{1}{\ref{c:beta:short}\cdot  g\cdot \sys_h(X)(g+|\log \sys_h(X)|)}.
	\end{eqnarray*}
	This completes the proof.
\end{proof}

As a direct consequence, we have
\begin{corollary}\label{cor:short:syst}
	Given a hyperbolic surface $X\in\T_g$ with homologically systolic length $\sys_h(X)\leq\left(\frac{\log g}{g}\right)^2$,  then 
$K(X)\geq \frac{\ref{a:shorter:systole}}{(\log g)^2}$ for some universal positive constant $\consta\label{a:shorter:systole}$.
\end{corollary}

\begin{proof}[Proof of Theorem \ref{thm:minimum}]
	Let $X_n\in\M_g$ be a sequence of hyperbolic surfaces realizing the infimum of the algebraic intersection form, that is 
	\begin{equation}\label{eq:limit:int}
		\lim_{n\to\infty} K(X_n)=\inf _{X\in\M_g} K(X).
	\end{equation}
	If $\{X_n\}_n$ is contained in some compact subset of $\M_g$, then the theorem follows by taking any accumulation point of $\{X_n\}$.  We now assume that $X_n\to\infty$ as $n\to\infty$.  By Lemma \ref{lem:nonsys}, we see that the length of nonseparating simple closed geodesics are bounded from below, say at least $\constc\label{c:min:sys:lower}>0$. That $X_n$ goes to infinity implies that the length of some separating simple closed geodesic goes to zero. 	Let $\Gamma_n$ be the set of separating simple closed geodesics on $X_n$ of length strictly less $\ref{c:min:sys:lower}$.  By Lemma \ref{lem:decreasing}, there exists a hyperbolic surface $X_n'$ such that
	\begin{itemize}
		\item $\ell_{X_n'}(\gamma)=\ref{c:min:sys:lower}$ for every $\gamma\in\Gamma_n$.
		\item  $\ell_{X_n'}(\alpha)\geq \ell_{X_n}(\alpha)\geq \ref{c:min:sys:lower}$ for every simple closed geodesic  $\alpha$ disjoint from any $\gamma\in \Gamma_n$.
	\end{itemize}
	Every simple closed geodesic on $X'_n$ that intersects some curve $\gamma\in \Gamma_n$ is of length at least twice the width of the collar neighborhood $\C(\gamma)$ of $\gamma$.
	Therefore, the sequence $\{X_n'\}$ is contained in a compact subset of $\M_g$.   Furthermore, by Lemma \ref{lem:near:bnd}, it follows that  $K(X_n')\leq K(X_n)$. Hence,  by \eqref{eq:limit:int}, the limit of any convergent subsequences of $\{X_n'\}$ attains the minimum $\min\limits_{X\in\M_g}K(X)$. This completes the proof.
\end{proof}

\section{Estimate of minimal intersection form}

In this section, we shall estimate the minimum of algebraic intersection form in the moduli space, proving Theorem \ref{thm:min:estimate}.  Before that, let us recall the following result by Balacheff, Parlier, and Sabourau.

\begin{theorem}[\cite{BPS2012}, Theorem 1.1]\label{thm:BPS}
	Let $M$ be a Riemannian surface of genus $g$ and area $4\pi(g-1)$. Then there exist $2g$ closed curves which induce a basis in $H_1(M,\Z)$ such that
	\begin{equation*}
		\ell_M(\alpha_k)\leq \frac{2^{16}}{\min\{\sys_h(M),1\}}\cdot  \frac{\log (2g-k+2)}{2g-k+1}\cdot g,
	\end{equation*}
	where $\sys_h(M)$ represents the homologically systolic length of $M$.
\end{theorem}
As a direct consequence, we have
\begin{proposition}\label{prop:large:systole}
	Let $M$ be a Riemannian surface of genus $g$ and area $4\pi(g-1)$. Then 
	\begin{equation*}
		K(M)\geq2^{-34} \left(\frac{\min\{\sys_h(M),1\}}{\log g}\right)^2,
	\end{equation*}
		where $\sys_h(M)$ represents the homologically systolic length of $M$.
\end{proposition}
\begin{proof}
	Let $\alpha_1,\cdots,\alpha_{2g}$ be the curves obtained from Theorem \ref{thm:BPS}. Since they induce a basis in $H_1(X,\Z)$, there exist two curves  $\alpha_{i_1}$ and $\alpha_{i_2}$ with $1\leq i_1\leq i_2\leq g+1$ such that $|\Int(\alpha_{i_1},\alpha_{i_2})|=1$.  Hence,
	\begin{eqnarray*}
		K(X)&\geq& \frac{|\Int(\alpha_{i_1},\alpha_{i_2})|}{\ell_X(\alpha_{i_1})\ell_X(\alpha_{i_2})}\geq \frac{1}{\ell_X(\alpha_{g+1})\ell_X(\alpha_{g+1})}\\
		&\geq&\left(\frac{\min\{\sys_h(X),1\}}{2^{16}}\right)^2 \frac{1}{\log ^2(g+1)}\\
		&\geq&\left(\frac{\min\{\sys_h(X),1\}}{2^{16}}\right)^2 \frac{1}{4\log ^2 g}=\frac{1}{2^{34}}\left(\frac{\min\{\sys_h(X),1\}}{\log g}\right)^2 .	\end{eqnarray*}
\end{proof}


\begin{proof}[Proof of Theorem \ref{thm:min:estimate}]
	\underline{\textbf{Upper bound}}. For any $X\in
	\T_g$,  let $\sys(X)$ (resp. $\sys_h(X)$) be the length of any (resp. homological) systole of $X$.  By \cite[Theorem 1.5]{MM2014}, we have
	\begin{equation*}
		K(X)\leq \frac{9}{\sys_h(X)^2}\leq \frac{9}{\sys(X)^2}.
	\end{equation*}
 Hence,
	\begin{equation*}
		\min_{X\in\M_g} K(X)\leq \frac{9}{\max_{X\in\M_g}\sys(X)^2}.
	\end{equation*}
	Recall that $\max_{X\in\M_g}\sys(X) \geq C\log g$ for some universal constant $C>0$ (\cite[Section 4]{BuserSarnak1994}, see also \cite{Brooks1999,KSV2007}; see \cite[Theorem 4.6]{Parlier2012} for a summary). Consequently, we have
	\begin{equation*}\label{eq:min:upper}
		\min_{X\in\M_g} K(X)\leq \frac{C}{(\log g)^2}.
	\end{equation*}
	 for some universal constant $C>0$.

	\underline{\textbf{Lower bound}}. Let $X\in \T_g$ be hyperbolic surface realizing the minimum: 
	\begin{equation*}
		K(X)=\min_{Y\in\M_g}K(Y).
	\end{equation*}	
	 If $\sys_h(X)\geq \ref{a:initial}$ or $\sys_h(X)<\left(\frac{\log g}{g}\right)^2$, then by Proposition \ref{prop:large:systole} or Corollary \ref{cor:short:syst}, we see that $K(X)\geq C (\log g)^{-2}$ for some universal constant $C$.  In the following, we assume that 
	 \begin{equation}\label{eq:sysh:assumption}
	 	 \left(\frac{\log g}{g}\right)^2\leq \sys_h(X)\leq \ref{a:initial}.
	 \end{equation}
	 	By Lemma \ref{lem:decreasing} and Lemma \ref{lem:near:bnd}, we may assume further that every separating simple closed geodesic on $X$ has length at least $\ref{a:initial}$. In particular, 
	 		 \begin{equation}\label{eq:sysh:assumption2}
	 	 \left(\frac{\log g}{g}\right)^2< \sys(X)< \ref{a:initial}.	 \end{equation} 	
	 	Let $\{\gamma_1,\cdots,\gamma_k\}$ be the set of simple closed geodesics on $X$ whose length are strictly less than $\ref{a:initial}$:
	 	\begin{equation}\label{eq:min:short:gamma}
\ell_X(\gamma_i)<\ref{a:initial}.
	 	\end{equation}
	 	In particular,  each $\gamma_i$ is nonseparating.

		Let $Y$ be the surface obtained in Lemma \ref{lem:auxiliary:surface}. Then $\sys_h(Y)\geq \sys(Y)=\ref{a:initial}$.  Combined with Proposition \ref{prop:large:systole}, this implies that
		\begin{equation*}
			K(Y)\geq \frac{2^{-34}\ref{a:initial}^2}{
			(\log g)^2}.
		\end{equation*}
	By Theorem \ref{thm:BPS},  there exist $2g$ oriented closed geodesics $\alpha_1,,\cdots,\alpha_{2g}$ on $Y$ which induces a basis of $H_1(Y,\Z)$ such that 
	\begin{equation}\label{eq:short:alpha:j}
		\ell_{Y}(\alpha_j)\leq \frac{2^{17}}{\ref{a:initial}}\cdot \log g, \quad\text{ for }1\leq j\leq g+1.
	\end{equation}
	  Being a basis of $H_1(Y,\Z)$ implies that two of the first $g+1$ curves, say $\alpha_1$ and $\alpha_2$, satisfy $|\Int(\alpha_1,\alpha_2)|= 1$.  Considering the intersection between $\alpha_i$ and $\gamma$, there are two cases.
      
      \textbf{Case 1:  $\sum_{j}(|\Int(\alpha_1,\gamma_j)|+|\Int(\alpha_2,\gamma_j)|)>0$.}  In this case, there exists one pair $(\alpha_i,\gamma_j)$ with positive algebraic intersection,  say $(\alpha_1,\gamma_1)$.  Then by Lemma \ref{lem:auxiliary:surface}, it follows that
      \begin{eqnarray*}
      	K(X)&\geq& \frac{|\Int(\alpha_1,\gamma_1)|}{\ell_X(\alpha_1)\ell_X(\gamma_1)}\\&\geq &
      	\frac{1}{\ref{a:fn}|\log \sys(X)|\ell_{Y}(\alpha_1)\ell_X(\gamma_i)}
      	\qquad(\text{ by Lemma }\ref{lem:auxiliary:surface} )\\
           	&\geq &  \frac{1}{\ref{a:fn}2^{18}(\log g)^2} 	\qquad(\text{ by } \eqref{eq:sysh:assumption2}, \eqref{eq:min:short:gamma}, \text{and } \eqref{eq:short:alpha:j}).
      \end{eqnarray*}

    \textbf{Case 2:  $\sum_{j}(|\Int(\alpha_1,\gamma_j)|+|\Int(\alpha_2,\gamma_j)|)=0$.}
	For each $\alpha_i$, applying Lemma \ref{lem:intersection:direction} for the each of the pair $(\alpha_i,\gamma_1)$ yields a set $A_{i}^1$ of oriented simple closed geodesics such that 
	\begin{itemize}
		\item $\cup_{\alpha'\in A_{i}^1} \alpha'$ is homologous to  $\alpha_i$,
		\item  $\sum_{\alpha'\in A_{i}^1} \ell_{Y}(\alpha')\leq \ell_Y(\alpha_i)$,
		\item $\cup_{\alpha'\in A_{i}^1} \alpha'$ is disjoint from $\gamma_1$.
	\end{itemize}
	Now inductively, suppose that for $1\leq j\leq k-1$, we have constructed the set $A_{i}^j$ satisfying the following:
	\begin{itemize}
		\item $\cup_{\alpha'\in A_{i}^j} \alpha'$ is homologous to  $\alpha_i$,
		\item  $\sum_{\alpha'\in A_{i}^j} \ell_{Y}(\alpha')\leq \ell_Y(\alpha_i)$,
		\item $\cup_{\alpha'\in A_{i}^j} \alpha'$ is disjoint from any of $\gamma_1,\cdots,\gamma_j$.
	\end{itemize}

	 For each $\alpha'\in A_i^j$, we  apply Lemma \ref{lem:intersection:direction} for  the pair $(\alpha',\gamma_{j+1})$, yielding a set $A_i^{j+1}$ of oriented simple closed geodesics such that
		\begin{itemize}
		\item $\cup_{\alpha''\in A_{i}^{j+1}} \alpha''$  is homologous to $\cup_{\alpha'\in A_{i}^j} \alpha'$, which is in turn homologous to   $\alpha_i$ by the inductive assumption,
		\item  $\sum_{\alpha''\in A_{i}^{j+1}} \ell_{Y}(\alpha'')\leq \sum_{\alpha'\in A_{i}^j} \ell_{Y}(\alpha')$, where the latter is at most $ \ell_Y(\alpha_i)$ by the inductive assumption,
		\item $\cup_{\alpha'\in A_{i}^{j+1}} \alpha'$ is disjoint from any of $\gamma_1,\cdots\gamma_j,\gamma_{j+1}$.
	\end{itemize}
After $k$ steps, we arrive at two sets $A_1^k$ and $A_2^k$ of oriented simple closed geodesics such that for each $i$, 
	\begin{itemize}
		\item $\cup_{\alpha'\in A_i} \alpha'$ is homologous to $\alpha_i$,
		\item $\sum_{\alpha'\in A_i}\ell_Y(\alpha')\leq \ell_Y(\alpha_i)$, and
		\item $\cup_{\alpha'\in A_i}\alpha'$ is disjoint from any of  $\gamma_1,\cdots,\gamma_k$.
	\end{itemize}
	Combining with the fact that $|\Int(\alpha_1,\alpha_2)|=1$, we infer that there exists $\alpha_i'\in A_i^k$ such that
	\begin{itemize}
		\item $|\Int(\alpha_1',\alpha_2')|>0$,
		\item $\ell_Y(\alpha_1')\leq \ell_Y(\alpha_1)$ and $\ell_Y(\alpha_2')\leq \ell_Y(\alpha_2)$,
		\item $\alpha_1'\cup\alpha_2'$ is disjoint from any of  $\gamma_1,\cdots,\gamma_k$.
	\end{itemize}
	Combined with Lemma \ref{lem:auxiliary:surface} and \eqref{eq:short:alpha:j}, this implies that
	\begin{equation*}
		\ell_X(\alpha_i')\leq \ell_Y(\alpha_i')\leq \frac{2^{17}}{\ref{a:initial}}\log g.
	\end{equation*}
Hence 
\begin{equation*}
	K(X)\geq \frac{|\Int(\alpha_1',\alpha_2')|}{\ell_X(\alpha_1')\ell_X(\alpha_2')}\geq \frac{2^{-34}\ref{a:initial}^2}{(\log g)^2}.
\end{equation*}
This completes the proof.
	\end{proof}

\section{Asymptotic behavior in moduli space}\label{sec:asymp}
In this section, we will describe the asymptotic behavior of the algebraic intersection form  $K(X)$ as the homologically systolic length  $\sys_h(X)$ goes to zero.  The key is to find \enquote{good representatives} of oriented simple closed curves approximating $K(X)$.
\subsection{Finding good representatives of curves}
We start with the following lemma.
\begin{lemma}\label{lem:intersection:11}
			Let $X\in\T_g$ be a hyperbolic surface with $\gamma_1,\cdots,\gamma_k$ being the set of nonseparating simple closed geodesics of length at most $\ref{a:initial}$. For any oriented simple closed geodesic $\alpha$ on $X$, there exists a set $A$ of oriented simple  closed geodesics  such that 
			\begin{enumerate}[(i)]
				\item  $\cup_{\alpha'\in A}~\alpha'$ is homologous to  $\alpha$, and
				\item $\i(\alpha',\gamma_j)\in\{0,1\}$ for every $\alpha'\in A$ and every $\gamma_j$,
				\item $\sum_{\alpha'\in A}\ell_X(\alpha') \leq \ell_X(\alpha) (1+\max\limits_{i}\frac{\ell_X(\gamma_i)}{|\log \ell_X(\gamma_i) |})$.
			\end{enumerate}
		\end{lemma}
		
		\begin{proof}
			Applying Lemma \ref{lem:intersection:direction} for every $\gamma_i$, we see that there exist oriented simple  closed geodesics $\alpha_1,\cdots,\alpha_m$ such that 
			\begin{enumerate}[(i)]
				\item  $\alpha_1\cup\alpha_2\cup\cdots\cup\alpha_m$ is homologous to  $\alpha$, and
				\item $\ell_X(\alpha_1)+\ell_X(\alpha_2)+\cdots+\ell_X(\alpha_m)\leq \ell_X(\alpha) $, and
				\item $|\Int(\gamma_j,\alpha_i)|=\i(\gamma,\alpha_i)$ for every $\alpha_i$ and every $\gamma_j$.
			\end{enumerate}
			In particular, the third property above means that for each pair $(i,j)$, $\alpha_i$ and $\gamma_j$ have the same algebraic intersection sign at all of their intersection points. Next, for each $\alpha_i$, we shall construct oriented simple closed geodesics with required properties.

			Let us start with $\alpha_1$.   Suppose that $\i(\alpha_1,\gamma_1)\geq2$. Otherwise, we move on to consider $\alpha_2$.  Let $\alpha_1:[0,\ell_X(\alpha_1)]\to X$ be a length parametrization of $\alpha_1$. Choose an orientation of $\gamma_1$ such that $\alpha_1$ enters the  collar neighbourhood $\C(\gamma_1)$ of $\gamma_1$ from the left side (see Section \ref{sec:thick:thin} for the definition of $\C(\gamma_1)$).  This is possible due to the property (iii) above. Furthermore, we assume that $\alpha_1(0)$ is on the left boundary component of the collar neighbourhood of $\gamma_1$, with repect to the orientation of $\gamma_1$.  Set $t_1=0$. Let $s_1\in (0,\ell_X(\alpha_1))$ be the first moment that $\alpha_1$ leaves   the  collar neighbourhood of $\gamma_1$. Let $t_2 \in (s_1,\ell_X(\alpha_1))$ be first moment that $\alpha_1$ enters the  collar neighbourhood of $\gamma_1$ again. Let $s_2\in (t_2,\ell_X(\alpha_1))$ be the first moment that $\alpha_1$ leaves the  collar neighborhood of $\gamma_1$ after $t_2$.  Let $c_i\in(t_i,s_i)$ be the moment that $\alpha_1|_{[t_i,s_i]}$ intersects $\gamma_1$.  Let $\gamma_1:[0,\ell_X(\gamma_1)]\to X$ be a length parametrization of $\gamma_1$. For $s,t\in[0,\ell_X(\gamma))$, let $\gamma_1|_{[s,t]}$ be the oriented segment from $\gamma_1(s)$ to $\gamma_1(t)$ with respect to the orientation of $\gamma_1$.  Let $\overline{\gamma_1}|_{[s,t]}$	 be the underlying segment of $\gamma_1|_{[s,t]}$	  with the \emph{opposite} orientation.	
			 Let $d_i\in [0,\ell_X(\gamma_1))$ be the moment that $\gamma_1$ intersects the segment $\alpha_1|_{[t_i,s_i]}$.   Let $\hat \alpha_{11}$ and $\hat \alpha_{12}$ be respectively the oriented simple closed curves defined as: 
			\begin{equation*}
				\hat \alpha_{11}:= \alpha_1|_{[t_1, c_1]}*\gamma_1|_{[d_1,d_2]}*\alpha_1|_{[c_2,t_1]} 
			\end{equation*}
			and
			\begin{equation*}
				\hat \alpha_{12}:= \alpha_1|_{[t_2, c_2]}*\overline{\gamma_1}|_{[d_1,d_2]}*\alpha_1|_{[c_1,t_2]} .
			\end{equation*}	
		  Let $\alpha'_{11}$ and $\alpha'_{12}$ be respectively the oriented simple closed geodesic homotopic to $\hat\alpha_{11}$ and $\hat\alpha_{12}$. Then
			 \begin{itemize}
				\item $\alpha'_{11}\cup\alpha'_{12}$ is homologous to $\alpha_1$,
				\item $1\leq \i(\alpha'_{11},\gamma_1)\leq \i(\alpha_1,\gamma_1)-1$,  $1\leq \i(\alpha'_{12},\gamma_1)\leq \i(\alpha_1,\gamma_1)-1$,
				\item $\ell_X(\alpha'_{11})+\ell_X(\alpha'_{12})\leq \ell_X(\alpha_1)+ 2\ell_X(\gamma_1)$,
			\end{itemize}
			where the last inequality follows from the fact that both ${\gamma_1}|_{[d_1,d_2]}$ and $\overline{\gamma_1}|_{[d_1,d_2]}$ are of length less than $\ell_X(\gamma_1)$.
			
			Repeating the above construction for each of $\alpha_{11}$ and $\alpha_{12}$ and the new resulting curves yields a set $A_{11}$ of oriented simple closed geodesics such that
			 \begin{itemize}
				\item $\cup_{\alpha'\in A_{11}}~\alpha'$ is homologous to $\alpha_1$,
				\item $ \i(\alpha',\gamma_1)=1$ every $\alpha'\in A_{11}$,
				\item $\sum_{\alpha'\in A_{11}}\ell_X(\alpha')\leq \ell_X(\alpha_1)+ 2\i(\alpha_1,\gamma_1)\ell_X(\gamma_1).$
			\end{itemize}
			
			Similarly, for every other $\alpha_i$ with $\i(\alpha_i,\gamma_1)\geq 1$, there exists a  set $A_{1i}$ of oriented simple closed geodesics such that
			 \begin{itemize}
				\item $\cup_{\alpha'\in A_{1i}}~\alpha'$ is homologous to $\alpha_i$,
				\item $ \i(\alpha',\gamma_1)=1$ every $\alpha'\in A_{1i}$,
				\item $\sum_{\alpha'\in A_{1i}}\ell_X(\alpha')\leq \ell_X(\alpha_i)+ 2\i(\alpha_i,\gamma_1)\ell_X(\gamma_1).$
			\end{itemize}
			Combining with the property (i), (ii), (iii) of $\alpha_1,\cdots,\alpha_m$ at the beginning of the proof, we arrive at a set $A_1=A_{11}\cup \cdots \cup A_{1m}$ of oriented simple closed geodesics such that
			 \begin{itemize}
				\item $\cup_{\alpha'\in A_{1}}~\alpha'$ is homologous to $\alpha$,
				\item $ \i(\alpha',\gamma_1)\in\{0,1\}$ every $\alpha'\in A_{1}$,
				\item $\sum_{\alpha'\in A_{1}}\ell_X(\alpha')\leq \ell_X(\alpha)+ 2\i(\alpha,\gamma_1)\ell_X(\gamma_1).$
			\end{itemize}
	Here the second item has an additional possibility $\i(\alpha',\gamma_1)=0$ is due to the fact that some of $\alpha_1,\cdots,\alpha_m$ may be disjoint from $\gamma_1$.

			Applying the above construction for  $\gamma_2, \gamma_3,\cdots,\gamma_k$ finitely many times, we arrive at a  set $A$ of oriented simple closed geodesics such that
			\begin{itemize}
				\item $\cup_{\alpha'\in A}~\alpha'$ is homologous to $\alpha$,
				\item $ \i(\alpha',\gamma_j)\in\{0,1\}$ for every $\alpha'\in A$ and every $\gamma_j$,
				\item $\sum_{\alpha'\in A}\ell_X(\alpha')\leq \ell_X(\alpha)+ 2\sum_j\i(\alpha,\gamma_j) \ell_X(\gamma_j) $.
			\end{itemize}
      The assumption $\ell_X(\gamma_j)\leq \ref{a:initial}$ implies that $\gamma_j$ admits  a collar neighbourhood of width 
	\begin{equation*}
		\arcsinh\frac{1}{\sinh (\ell_X(\gamma_j)/2)}\geq |\log\ell_X(\gamma_j)|.
	\end{equation*}
     Hence, the length $\ell_X(\alpha)$ satifies 
       \begin{equation*}
      	\ell_X(\alpha)\geq \sum_j 2\cdot \i(\alpha,\gamma_j)\cdot |\log \ell_X(\gamma_j)|.
      \end{equation*}
      Therefore,
      \begin{eqnarray*}
      	2\sum_j\i(\alpha,\gamma_j) \ell_X(\gamma_j)&\leq&  \max_j \frac{\ell_X(\gamma_j)}{|\log \ell_X(\gamma_j)|}\sum_j 2\cdot \i(\alpha,\gamma_j)\cdot |\log \ell_X(\gamma_j)|\\
      &	\leq&  \max_j \frac{\ell_X(\gamma_j)}{|\log \ell_X(\gamma_j)|}\cdot \ell_X(\alpha).
      \end{eqnarray*}
      Accordingly, 
   \begin{equation*}
   \sum_{\alpha'\in A}\ell_X(\alpha')\leq \left(1+  \max_j \frac{\ell_X(\gamma_j)}{|\log \ell_X(\gamma_j)|}\right)\cdot \ell_X(\alpha).
   \end{equation*}
   This completes the proof.
		\end{proof}

		Next, we shall compare curves that wrap many times about short curves with those that wrap a few times. 	As before, let  $\gamma_1,\cdots,\gamma_k$ be the set of nonseparating simple closed geodesics of length at most   $\ref{a:initial}$ on $X\in\T_g$.  Recall that the collar neighbourhood  $\C(\gamma_i)$ of every short curve $\gamma_i$ is foliated by its \emph{perpendiculars}, geodesic arcs connecting one boundary component to the other and intersects $\gamma_i$ perpendicularly (see Section \ref{sec:thick:thin}). 		Let $\alpha_0$ be an oriented simple closed geodesic on $X$ which crosses $\gamma_i$ at most once and which crosses every perpendicular of $\C(\gamma_i)$ at most once.    Let $(m_i)_{1\leq i\leq k}$ be a $k$-tuple of integers such that $m_i\neq0$ if and only if $\i(\alpha_0,\gamma_i)=1$.  Let $\alpha$ be a composition of $m_j$-th Dehn twist $T^{m_j}_{\gamma_j}$ of $\alpha_0$ around  $\gamma_j$: 
		\begin{equation*}
			\alpha=T^{m_k}_{\gamma_k}\circ\cdots\circ T^{m_1}_{\gamma_1}(\alpha_0).
		\end{equation*} 
		\begin{lemma}\label{lem:small:twist:11}
		Let $X$, $\gamma_1,\cdots,\gamma_k$, $\alpha_0$, and $\alpha$ be as above.  Then, after choosing appropriate orientation for each $\gamma_i$, we have
			\begin{enumerate}[(i)]
				\item  $\alpha_0\cup m_1\gamma_1\cup\cdots\cup m_k\gamma_k$ is homologous to  $\alpha$, and
				\item $\ell_X(\alpha_0)+|m_1|\ell_X(\gamma_1)+\cdots+|m_k|\ell_X(\gamma_k)\leq \ell_X(\alpha) (1+\max\limits_i \frac{\ref{a:small:twists}}{|\log \ell_X(\gamma_i)|})$,\end{enumerate}
				for some universal constant $\consta\label{a:small:twists}$.
		\end{lemma}
		
		\begin{proof}
		The first item is straightforward. For the second item, consider the collar neighborhood $\C(\gamma_i)$ of $\gamma_i$. 
		Note that by assumption, $\alpha$ crosses $\gamma_i$ at most once. Each of $\alpha_0\cap \C(\gamma_i)$ and $\alpha\cap \C(\gamma_i)$  is either an empty set or a single segment. Moreover, the segment (if any) $\alpha_0\cap \C(\gamma_i)$ crosses every perpendicular of $\C(\gamma_i)$ at most once, and hence $\alpha\cap \C(\gamma_i)$ crosses every perpendicular of $\C(\gamma_i)$ at most $|m_i|+1$ times. 
		Then,  by \eqref{eq:thin:length}, we see that 
		\begin{equation*}
		\left|	\ell_X(\alpha\cap \C(\gamma_i) )- (\ell_X(\alpha_0\cap \C(\gamma_i)+ |m_i|\ell_X(\gamma_i))\right|\leq C 
		\end{equation*}
		for some universal constant $C$. Let $X_{\thick}:=X\backslash(\C(\gamma_1)\cup\cdots\cup\C(\gamma_k))$ be the thick part of $X$. Note that each of $\alpha_0\cap X_\thick $ and $\alpha\cap X_\thick $  has $\sum_i\i(\alpha,\gamma_i)$ segments. Furthermore, they are homotopic to each other relative to the boundary components of $\cup_i\C(\gamma_i)$. Since $\ell_X(\gamma_i)\leq \ref{a:initial}$, it follows that each boundary component of $\cup_i\C(\gamma_i)$ is of length most $2\sqrt{2}$ (see Section \ref{sec:thick:thin}). Therefore, 
		\begin{equation*}
			\left|\ell_X(\alpha_0\cap X_\thick)- \ell_X(\alpha\cap X_\thick)\right|\leq 2\sqrt{2}\sum_i \i(\alpha,\gamma_i).
		\end{equation*}
		On the other hand, the assumption $\gamma_i\leq \ref{a:initial}$ implies that the width of the collar neighbourhood $\C(\gamma_i)$ satisfies
	\begin{equation*}
		\arcsinh\frac{1}{\sinh (\ell_X(\gamma_i)/2)}\geq |\log\ell_X(\gamma_i)|.
	\end{equation*}
Hence, 
		\begin{equation*}
			\ell_X(\alpha)\geq \sum_i\ell_X(\alpha\cap \C(\gamma_i))\geq 2\sum_i \i(\alpha,\gamma_i)|\log \ell_X(\gamma_i)|.
		\end{equation*}
		Combining the estimates above yields
		\begin{eqnarray*}
			\ell_X(\alpha_0)+\sum_i |m_i|\ell_X(\gamma_i)&\leq& \ell_X(\alpha)+\sum_i \i(\alpha,\gamma_i)(2\sqrt{2}+C)\\
			&\leq & \ell_X(\alpha) \left(1+\max_{i}\frac{2\sqrt{2}+C}{|\log \ell_X(\gamma_i)|}\right).
		\end{eqnarray*}
		This completes the proof.
		\end{proof}
		
		\begin{proof}[Proof of Theorem \ref{thm:asymp}] 
		 For any constant $\epsilon>0$, let $\constc\label{c:small:epsilon}$ be a constant depending only $\epsilon$ such that  
		 \begin{equation}\label{eq:asymp:delta}
		 	\ref{c:small:epsilon}<1,~\frac{1}{|\log\ref{c:small:epsilon}|}<\epsilon, ~\frac{\ref{a:small:twists}}{|\log\ref{c:small:epsilon}|}<\epsilon.
		 \end{equation}
		   Let $X\in\M_g$ be a hyperbolic surface with $\sys_h(X)<  \ref{c:small:epsilon}$. Let $\alpha$ and $\beta$ be a pair of oriented simple closed geodesics on $X$ such that 
		 \begin{equation}\label{eq:asymp:K:ab}
		 	(1-\epsilon) \cdot K(X)\leq \frac{|\Int(\alpha,\beta)|}{\ell_X(\alpha)\ell_X(\beta)}.
		 \end{equation}
		 Let $\Gamma:=\{\gamma_1,\cdots,\gamma_k\}$ be the set of nonseparating simple closed geodesics on $X$ of length at most $\ref{c:small:epsilon}$:
		 \begin{equation}\label{eq:aysmp:short}
		 	\ell_X(\gamma_j)\leq \ref{c:small:epsilon}.
		 \end{equation}
			Let $A$ (resp.  $B$) be the set of oriented simple closed geodesics obtained by applying Lemma \ref{lem:intersection:11} to $\alpha$ (resp. $\beta$).  Then
			\begin{eqnarray}
			\nonumber	\frac{|\Int(\alpha,\beta)|}{\ell_X(\alpha)\ell_X(\beta)}&\leq & \frac{\left|\sum\limits_{\alpha'\in A,~\beta'\in B}\Int(\alpha',\beta')\right|}{\left(\sum\limits_{\alpha'\in A}\ell_X(\alpha')\right)\left(\sum\limits_{\beta'\in A}\ell_X(\beta')\right)}\cdot (1+\max_j\frac{\ell_X(\gamma_j)}{|\log \ell_X(\gamma_i)|})^2
			\\
				&\leq & \max\limits_{\alpha'\in A,~\beta'\in B} \frac{|\Int(\alpha',\beta')|}{\ell_X(\alpha')\ell_X(\beta')}\cdot (1+\epsilon) ^2\label{eq:asymp:ab} .\qquad(\text{ by } \eqref{eq:asymp:delta}, \eqref{eq:aysmp:short})
			\end{eqnarray}
			Note that 
						 \begin{equation}
			 	\label{eq:asymp:int:abg}
			 	|\Int(\alpha',\gamma_i)|=\i(\alpha',\gamma_i)\leq1,\quad 			 	|\Int(\beta',\gamma_i)|=\i(\beta',\gamma_i)\leq 1.
			 \end{equation}
			 For each $\gamma_i\in\Gamma$, let $T^n_{\gamma_i}$ be the $n$-th Dehn twist around $\gamma_i$. Let $\alpha''$ be one of the shortest simple closed geodesics in the Dehn twists orbit
			\begin{equation*}
			\{	T^{n_1}_{\gamma_1}\circ T^{n_2}_{\gamma_2} \circ\cdots\circ T^{n_k}_{\gamma_k} (\alpha'):n_i\in \mathbb Z\}.
			\end{equation*}	
			Assume that $\alpha''=	T^{n_1}_{\gamma_1}\circ \cdots\circ T^{n_2}_{\gamma_2} \circ T^{n_k}_{\gamma_k} (\alpha')$.
			Similarly define $\beta'':=T^{m_1}_{\gamma_1}\circ T^{m_2}_{\gamma_2} \circ \cdots\circ T^{m_k}_{\gamma_k} (\beta')$. Being shortest in their Dehn twists orbits implies that both $\alpha''$ and $\beta''$ intersect every perpendicular of the collar neighbourhood $\C(\gamma_i)$ at most once.  Applying  Lemma \ref{lem:small:twist:11} to $(\alpha'',\alpha')$ and $(\beta'',\beta')$ independently, we see that
			\begin{eqnarray*}
			\nonumber&&\frac{|\Int(\alpha',\beta')|}{\ell_X(\alpha')\ell_X(\beta')}(1+\epsilon)^{-2}\\
				\nonumber&\leq & \frac{|\Int(\alpha',\beta')|}{\ell_X(\alpha')\ell_X(\beta')}\cdot(1+\max_j \frac{\ref{a:small:twists}}{|\log \ell_X(\gamma_j)|})^{-2} \qquad(\text{ by } \eqref{eq:asymp:delta}, \eqref{eq:aysmp:short})
				\\
				\nonumber&\leq &\frac{|\Int(\alpha''\cup (\cup_i n_i\gamma_i)),\beta''\cup(\cup_j m_j\gamma_j))|}{(\ell_X(\alpha'')+\sum_i|n_i|\ell_X(\gamma_i))(\ell_X(\beta'')+\sum_j|m_j|\ell_X(\gamma_j))}
				\\
				\nonumber&\leq & \frac{\sum_j |m_j| \cdot |\Int(\alpha'',\gamma_j)|+\sum_i|m_i|\cdot |\Int(\beta'',\gamma_i))|+|\Int(\alpha'',\beta'')|}{\sum_j |m_j|\ell_X(\alpha'')\ell_X(\gamma_j) +\sum_i |n_i|\ell_X(\beta'')\ell_X(\gamma_i)+\ell_X(\alpha'')\ell_X(\beta'')}
\\ \nonumber&\leq & \max_{i,j}\left\{\frac{|\Int(\alpha'',\gamma_j)|}{\ell_X(\alpha'',\gamma_j)},\frac{|\Int(\beta'',\gamma_i)|}{\ell_X(\beta'')\ell_X(\gamma_i)}, \frac{|\Int(\alpha'',\beta'')|}{\ell_X(\alpha'')\ell_X(\beta'')}\right\}\\
\nonumber&\leq & \max \left\{ \max_{i}\frac{1}{\ell_X(\alpha_{\gamma_i})\ell_X(\gamma_i)}, \frac{|\Int(\alpha'',\beta'')|}{\ell_X(\alpha'')\ell_X(\beta'')}\right\},
			\end{eqnarray*}	
			where the last inequality follows from \eqref{eq:asymp:int:abg} and   the assumption that $\alpha_{\gamma_i}$ is a shortest nonseparating simple closed geodesic which intersects $\gamma_i$ exactly once.

			Next, we shall compare 
			\begin{equation*}
				\max_i \frac{1}{\ell_X(\alpha_{\gamma_i})\ell_X(\gamma_i)} \quad\text{ and } \quad\frac{|\Int(\alpha'',\beta'')|}{\ell_X(\alpha'')\ell_X(\beta'')}.
			\end{equation*}
			To begin, suppose $\gamma_1$ is one of the shortest one among $\gamma_1,\cdots,\gamma_k$, that is, $\ell_X(\gamma_1)=\sys_h(X)$ . Since $\alpha_{\gamma_1}$ is one of the shortest simple closed geodesic intersecting $\gamma_1$ exactly once,  then by \eqref{eq:nonsys:beta}, we see that
				\begin{equation*}		\ell_X(\alpha_{\gamma_1})\leq\ref{c:beta:short}\cdot  g (g+|\log{\sys_h(X)}|)
	\end{equation*}	
	for some universal constant $\ref{c:beta:short}$.
	Hence 	
		\begin{equation*}
		\max_i	\frac{1}{\ell_X(\alpha_{\gamma_i})\ell_X(\gamma_i)} \geq 	\frac{1}{\ell_X(\alpha_{\gamma_1})\ell_X(\gamma_1)} \geq \frac{1}{\ref{c:beta:short}\cdot  g \cdot \sys_h(X)(g+|\log{\sys_h(X)}|)}
			\end{equation*}
			which goes to infinity as $\sys_h(X)\to0$. On the other hand,  we
			claim that 
			$				 \frac{|\Int(\alpha'',\beta'')|}{\ell_X(\alpha'')\ell_X(\beta'')}$  is bounded from above.
		 To see this, let $\C(\gamma_j)$ be the collar neighbourhood of $\gamma_j$. Let $X_\thin:=\cup_i \C(\gamma_i)$ and $X_\thick=X\backslash X_\thin$  
		 be respectively the thin part and the thick part of $X$ respectively. Then any simple closed curve in $X_\thick$ is of length bounded from below by some positive constant $\constc\label{c:bounded:below}$.  Applying an argument similar as the proof of \cite[Proposition 3.4]{MM2014} yields that
			\begin{equation*}
				\frac{|\Int(\alpha''\cap X_\thick,\beta''\cap X_\thick)|}{\ell_X(\alpha''\cap X_\thick)\ell_X(\beta''\cap X_\thick)}\leq \frac{9}{\ref{c:bounded:below}^2}.
			\end{equation*}
			Note that   $\i(\alpha'',\gamma_j)=\i(\alpha',\gamma_j)\leq 1$ and $\i(\beta'',\gamma_j)=\i(\beta',\gamma_j)\leq1$. Furthermore, both $\alpha''$ and $\beta''$ intersect each perpendicular of the collar neighbourhood $\C(\gamma_j)$ at most once. Therefore, $\alpha''\cap\C(\gamma_j)$ intersects $\beta''\cap \C(\gamma_j)$ at most twice:
			\begin{equation*}
				\i(\alpha''\cap\C(\gamma_j),\beta''\cap\C(\gamma_j))\leq 2\cdot  \i(\alpha'',\gamma_j)\i(\beta'',\gamma_j).
			\end{equation*}
			Note that the collar neighbourhood $\C(\gamma_i)$ of $\gamma_i$ is of width at least $|\log \ell_X(\gamma_i)|$. Hence,
			\begin{equation*}
				\ell_X(\alpha''\cap X_\thin)\geq \sum_j 2\i(\alpha'',\gamma_j)|\log \ell_X(\gamma_j)|
			\end{equation*}
			and 
						\begin{equation*}
				\ell_X(\beta''\cap X_\thin)\geq \sum_j 2\i(\beta'',\gamma_j)|\log \ell_X(\gamma_j)|.
			\end{equation*}

			 Accordingly, 
			\begin{eqnarray*}
				&&	\frac{|\Int(\alpha''\cap X_\thin,\beta''\cap X_\thin)|}{\ell_X(\alpha''\cap X_\thin)\ell_X(\beta''\cap X_\thin)}\\&\leq& \frac{\sum_j 2 \i(\alpha'',\gamma_j)\i(\beta'',\gamma_j)}{\ell_X(\alpha''\cap X_\thin)\ell_X(\beta''\cap X_\thin)}\\
					&\leq & \frac{\sum_j 2\i(\alpha'',\gamma_j)\i(\beta'',\gamma_j)}{\sum_j 2\i(\alpha'',\gamma_j)|\log \ell_X(\gamma_j)|\cdot \sum_j 2\i(\beta'',\gamma_j)|\log \ell_X(\gamma_j)|}\\
					&\leq &\max_{j}  \frac{ 2 }{ 2|\log \ell_X(\gamma_j)|\cdot 2|\log \ell_X(\gamma_j)|}\leq \frac{1}{2|\log \ref{c:small:epsilon}|^2}.
			\end{eqnarray*}
			Therefore,
			\begin{eqnarray*}
			&&	\frac{|\Int(\alpha'',\beta'')|}{\ell_X(\alpha'')\ell_X(\beta'')}\\ &\leq & 	\frac{|\Int(\alpha''\cap X_\thick,\beta''\cap X_\thick)|+|\Int(\alpha''\cap X_\thin,\beta''\cap X_\thin)|}{(\ell_X(\alpha''\cap X_\thick)+\ell_X(\alpha''\cap X_\thin ))(\ell_X(\beta''\cap X_\thick)+\ell_X(\beta''\cap X_\thin )}\\
			&\leq & \max\left\{\frac{|\Int(\alpha''\cap X_\thick,\beta''\cap X_\thick)|}{\ell_X(\alpha''\cap X_\thick)\ell_X(\beta''\cap X_\thick)}, 	\frac{|\Int(\alpha''\cap X_\thin,\beta''\cap X_\thin)|}{\ell_X(\alpha''\cap X_\thin)\ell_X(\beta''\cap X_\thin)}\right \}\\
			&\leq & \max\{ \frac{9}{\ref{c:bounded:below}^2},\frac{1}{2|\log \ref{c:small:epsilon}|^2}\}.
			\end{eqnarray*}
			This completes the proof of the claim. Consequently, there exists a constant $\constc<\ref{c:small:epsilon}\label{c:shorter:epsilon}$ such that if $\sys_h(X)\leq \ref{c:shorter:epsilon}$ then
			\begin{equation*}
				\frac{|\Int(\alpha',\beta')|}{\ell_X(\alpha')\ell_X(\beta')}\leq  \max_i \left\{ \frac{1}{\ell_X(\alpha_{\gamma_i})\ell_X(\gamma_i)}\right\}\cdot(1+\epsilon)^2.
			\end{equation*}		
		Combining this with \eqref{eq:asymp:K:ab} and \eqref{eq:asymp:ab}, we infer that
		\begin{equation*}
			K(X)\leq  \max_{i}\left\{\frac{1}{\ell_X(\alpha_{\gamma_i})\ell_X(\gamma_i)}\right\}\cdot\frac{(1+\epsilon)^4}{1-\epsilon}\leq \max_{\gamma\in\Gamma}\left\{\frac{1}{\ell_X(\alpha_{\gamma})\ell_X(\gamma)}\right\}\cdot\frac{(1+\epsilon)^4}{1-\epsilon},
		\end{equation*}	
		where, by assumption, $\Gamma$ is the set of simple closed geodesics of length at most 1 and $\alpha_\gamma$ is one of the shortest simple closed geodesics intersecting $\gamma\in\Gamma$ exactly once.  
			On the other hand, the definition of $K(X)$ implies that 
			\begin{equation*}
				K(X)\geq \max_{\gamma\in\Gamma}\left\{\frac{1}{\ell_X(\alpha_{\gamma})\ell_X(\gamma)}\right\}.
			\end{equation*}
			Hence, as $\sys_h(X)\to0$, we have 
			\begin{equation*}
			1\leq \frac{K(X)}{\max\limits_{\gamma\in\Gamma}\left\{\frac{1}{\ell_X(\alpha_{\gamma})\ell_X(\gamma)}\right\}}\leq \frac{(1+\epsilon)^4}{1-\epsilon}.
			\end{equation*}
			The theorem then follows from the arbitrariness of $\epsilon$. 
				\end{proof}

\subsection{Examples}\label{sec:exapmple}   

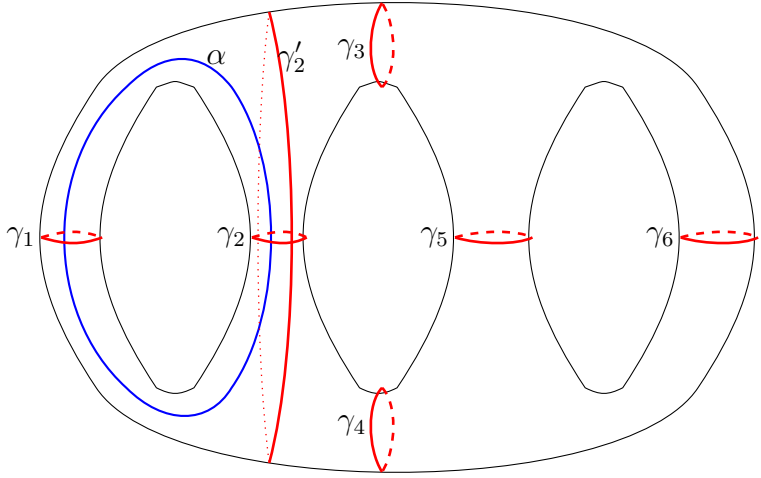
\begin{figure}
	\begin{tikzpicture}
		\draw (-4,2)..controls (-3,3.5) and (3,3.5).. (4,2).. controls (5,0.5) and (5,-0.5).. (4,-2)..controls (3,-3.5) and (-3,-3.5)..(-4,-2)..controls (-5,-0.5) and (-5,0.5)..(-4,2);

		\draw (4-1-2.5-0.5,2).. controls (5-1-2.5-0.5,0.5) and (5-1-2.5-0.5,-0.5).. (4-1-2.5-0.5,-2) ..controls (4-1-2.5-0.2-0.5,-2.1) and (4-1-2.5-0.3-0.5,-2.1).. (-4+6.5-2.5-0.5,-2)..controls (-5+6.5-2.5-0.5,-0.5) and (-5+6.5-2.5-0.5,0.5)..(-4+6.5-2.5-0.5,2)..controls ((4-1-2.5-0.2-0.5,2.1)and (4-1-2.5-0.3-0.5,2.1)..(4-1-2.5-0.5,2);
		 
		 	\draw (4-1,2).. controls (5-1,0.5) and (5-1,-0.5).. (4-1,-2) ..controls (4-1-0.2,-2.1) and (4-1-0.3,-2.1).. (-4+6.5,-2)..controls (-5+6.5,-0.5) and (-5+6.5,0.5)..(-4+6.5,2)..controls ((4-1-0.2,2.1)and (4-1-0.3,2.1)..(4-1,2);
		 	
		 		\draw (4-1-6+0.3,2).. controls (5-1-6+0.3,0.5) and (5-1-6+0.3,-0.5).. (4-1-6+0.3,-2) ..controls (4-1-6-0.2+0.3,-2.1) and (4-1-6-0.3+0.3,-2.1).. (-4+6.5-6+0.3,-2)..controls (-5+6.5+0.3-6,-0.5) and (-5+6.5+0.3-6,0.5)..(-4+6.5-6+0.3,2)..controls ((4-1-6-0.2+0.3,2.1)and (4-1-6-0.3+0.3,2.1)..(4-1-6+0.3,2);
		 		
		 		\draw[blue,line width=0.8pt] (-3.6, 2)..controls (-3.1,2.5) and (-2.6, 2.5)..(-2.2,2) ..controls (-1.5,1) and (-1.5,-1)..(-2.2,-2)..controls(-2.5,-2.5) and (-3.1,-2.5) .. (-3.6,-2)..controls (-4.7,-1) and (-4.7,1)..(-3.6,2); 
		 		
		 		\draw[red,line width=1pt](-4.73,0)..controls (-4.4,-0.1) and (-4.2,-0.1)..(-3.92,0);
		 		\draw[red,line width=1pt,dashed](-4.73,0)..controls (-4.4,0.1) and (-4.2,0.1)..(-3.92,0);
		 		
		 		\draw[red,line width=1pt](-4.73+2.8,0)..controls (-4.4+2.8,-0.1) and (-4.2+2.8,-0.1)..(-3.9+2.7,0);
		 		\draw[red,line width=1pt,dashed](-4.73+2.8,0)..controls (-4.4+2.8,0.1) and (-4.2+2.8,0.1)..(-3.9+2.7,0);
		 		
		 		\draw[red,line width=1pt](-4.73+2.8+2.7,0)..controls (-4.4+2.8+2.7,-0.1) and (-4.2+2.8+3,-0.1)..(-3.9+2.7+3,0);
		 		\draw[red,line width=1pt,dashed](-4.73+2.8+2.7,0)..controls (-4.4+2.8+2.7,0.1) and (-4.2+2.8+3,0.1)..(-3.9+2.7+3,0);
		 		
		 		\draw[red,line width=1pt](-4.73+2.8+2.7+3,0)..controls (-4.4+2.8+2.7+3,-0.1) and (-4.2+2.8+3+3,-0.1)..(-3.9+2.7+3+3,0);
		 		\draw[red,line width=1pt,dashed](-4.73+2.8+2.7+3,0)..controls (-4.4+2.8+2.7+3,0.1) and (-4.2+2.8+3+3,0.1)..(-3.9+2.7+3+3,0);
		 		\draw[red,line width=1pt] (0-0.2,3.12)..controls(-0.2-0.2,2.8) and (-0.2-0.2,2.2)..(0-0.2,2);
		 		\draw[red, line width=1pt,dashed](0-0.2,2)..controls (0.2-0.2,2.2) and (0.2-0.2,2.8).. (0-0.2,3.12);
		 		\draw[red,line width=1pt] (0-0.2,-3.12)..controls(-0.2-0.2,-2.8) and (-0.2-0.2,-2.2)..(0-0.2,-2);
		 		\draw[red, line width=1pt,dashed](0-0.2,-2)..controls (0.2-0.2,-2.2) and (0.2-0.2,-2.8).. (0-0.2,-3.12);
		 		
		 		\draw[dotted,red,line width=0.5pt] (-1.7,3)..controls (-1.9,1.5) and (-1.9,-1.5) ..(-1.7,-3);
		 		\draw[red,line width=1pt] (-1.7,3)..controls (-1.3,1.5) and (-1.3,-1.5) ..(-1.7,-3);
		 		
		 		\draw (-5,0)node{$\gamma_1$}
		 		(-2.2,0)node{$\gamma_2$}
		 		(0.5,0)node{$\gamma_5$}
		 		(3.5,0)node{$\gamma_6$}
		 		(-0.6,2.5)node{$\gamma_3$}
		 		(-0.6,-2.5)node{$\gamma_4$}
		 		(-2.4,2.4)node{$\alpha$}
		 		(-1.4,2.4)node{$\gamma_2'$};
	
		\end{tikzpicture}
		\caption{A pants decomposition $\{\gamma_1,\gamma_2,\cdots,\gamma_6\}$ consisting of nonseparating simple closed geodesics on a hyperbolic surface of genus three.}
		\label{fig:example}
\end{figure}
\begin{example}
Let $\gamma_1,\cdots,\gamma_{3g-3}$ be a collection of nonseparating simple closed geodesics as in Figure \ref{fig:example}. Let $0\leq\delta\leq1$ be a fixed constant. Let $X_n\in\M_g$ be a hyperbolic surface such that $\ell_{X_n}(\gamma_1)=\frac{1}{n^\delta}$ and $\ell_{X_n}(\gamma_{i})=\frac{1}{n}$ for $i\geq 2$.  As $n\to\infty$, the geodesics  $\gamma_2,\cdots,\gamma_{3g-3}$ become the homological systoles of $X_n$. The geodesic $\gamma_1$ has a collar neighbourhood of width $\delta\log n$ while $\gamma_i$ has a collar neighbourhood of width $ \log n$ for $i\geq2$. Any nonseparating simple closed geodesic intersecting distinct from any of $\{\gamma_i\}$ has to intersect at least two distinct curves of $\{\gamma_i\}$, and hence is of length at least $2(1+\delta)\log n=2(1+\delta)|\log \sys_h(X_n)|$. By  Theorem \ref{thm:asymp}, we see that 
		\begin{equation*}
			\lim_{n\to\infty} {K(X_n)}\cdot {2(1+\delta) \cdot \sys_h(X_n) |\log \sys_h(X_n)|}\leq 1.
		\end{equation*}
		On the other hand, let $\alpha$ be one of the shortest simple closed geodesics which intersect only $\gamma_1$ and $\gamma_2$ as in Figure \ref{fig:example}. Then  
		$\i(\alpha,\gamma_2)=1$ and  $|\ell_{X_n}(\alpha)-2(1+\delta)\log n|\leq C$ for some universal constant $C$ (by \eqref{eq:thin:length}). Therefore,
		\begin{eqnarray*}
			K(X_n)&\geq & \frac{|\Int(\alpha,\gamma_2)|}{\ell_{X_n}(\alpha)\ell_{X_n}(\gamma_2)}
			\\
			&\geq& \frac{1}{(2(1+\delta) \log n+C)\cdot \ell_{X_n}(\gamma_2)} \\
			&= &\frac{1}{(2(1+\delta) |\log \sys_h(X_n)|+C)\cdot \sys_h(X_n)},
		\end{eqnarray*}
		where in the last equation we use the assumption that $\ell_{X_n}(\gamma_2)=\sys_h(X_n)=\frac{1}{n}$.
		Letting $n\to\infty$, we see that 
		\begin{equation*}
			\lim_{n\to\infty} {K(X_n)}\cdot {2(1+\delta) \cdot \sys_h(X_n) |\log \sys_h(X_n)|}\geq  1.
		\end{equation*}
		In summary, we have
		\begin{equation*}
			\lim_{n\to\infty} {K(X_n)}\cdot {2(1+\delta) \cdot \sys_h(X_n) |\log \sys_h(X_n)|}= 1.
		\end{equation*}
		This example shows that the asymptotic behavior of the algebraic intersection form cannot be described by a single function of the homologically systolic length.
\end{example}

\begin{example}\label{example:disjoint:systole}
We continue with the notation in the previous example. Let $\gamma_2'$ be the simple closed geodesic as in Figure \ref{fig:example}, let $\ell_{X_n}(\gamma_1)=\frac{1}{n^{1-1/n}}$, $\ell_{X_n}(\gamma_2')=\frac{1}{n}$, $\ell_{X_n}(\gamma_i)=\frac{1}{n}$ for $i\geq 3$. Then, as $n\to\infty$, 
	\begin{equation*}
		\frac{1}{\ell_{X_n}(\alpha)\ell_{X_n}(\gamma_1)}\sim \frac{n^{1-1/n}}{2 (1-1/n) |\log n| }\sim \frac{n}{2|\log n|}.
	\end{equation*}
On the other hand, for each homological systole $\gamma_i$ ($i\geq 3$), any nonseparating simple closed geodesic $\alpha_{\gamma_i}$ intersecting $\gamma_i$  once has to intersect at least  two distinct curves among $\gamma_2',\gamma_3,\cdots,\gamma_{3g-3}$.  Hence $\ell_X(\alpha_{\gamma_j})\geq 4\log n$. Accordingly, 
\begin{equation*}
	\max_{3\leq i\leq 3g-3}\frac{1}{\ell_{X_n}(\alpha_{\gamma_i})\ell_{X_n}(\gamma_i)}\leq \frac{n}{4|\log n|}.
\end{equation*}
Therefore, by Theorem \ref{thm:asymp}, as $n\to\infty$, 
	\begin{equation*}
		K(X)\sim \frac{n}{2|\log n|}.
	\end{equation*}
This example shows that considering only homological systoles and curves intersecting them is not sufficient for describing the asymptotic behavior of the algebraic intersection form. 
\end{example}

\bibliographystyle{alpha}
\bibliography{Intersection.bib}

\end{document}